\documentclass[11pt]{amsart}
\usepackage{amssymb,amsthm}
\usepackage[curve]{xypic}

\begin{document}

\hyphenation{char-ac-ter-iza-tion}
\hyphenation{com-pact-i-fi-ca-tion}
\hyphenation{com-pact-i-fi-ca-tions}
\hyphenation{con-fig-u-ra-tion}
\hyphenation{cor-re-spond-ing}
\hyphenation{el-e-ment}
\hyphenation{el-e-ments}
\hyphenation{ex-cep-tion-al}
\hyphenation{in-ter-sec-tion}
\hyphenation{in-ter-sec-tions}
\hyphenation{pa-ram-e-trized}
\hyphenation{par-ti-tion}
\hyphenation{poly-di-ag-o-nal}
\hyphenation{poly-no-mi-al}
\hyphenation{si-mul-ta-neous}
\hyphenation{sub-lat-tice}

\newtheorem{thm}{Theorem}
\newtheorem{lem}{Lemma}		
\newtheorem{prop}{Proposition}
\newtheorem{corol}{Corollary}
\newtheorem*{FBL}{Flag Blowup Lemma}

\theoremstyle{remark}
\newtheorem*{rem}{Remark}
\newtheorem*{defn}{\bf Definition}
\newtheorem*{examples}{\bf Examples}

\renewcommand\theenumi{\alph{enumi}}
\renewcommand\labelenumi{(\theenumi)}



\def\inquotes#1{`#1'}
\def\indoublequotes#1{``#1''}

%
\def\ExtOn{
	\ifnum\catcode96=13 	
	\else
		\toks0={`}
		\catcode21=11	
		\catcode``=13
	\fi
	}
\ExtOn

%
\def`#1{
	\ifmmode\csname ^^U#1\endcsname
	\else\the\toks0#1\fi
	}

%
\def\defu#1#2{
	\edef\0{
		\def\csname ^^U#1\endcsname{#2}
		}
	\0
	}

%
\def\ExtOff{
	\catcode96=12	 
	\catcode21=12 
	}


\defu a{\alpha}
\defu b{\beta}
\defu c{\chi}
\defu d{\delta}
\defu e{\epsilon}
\defu f{\phi}
\defu g{\gamma}
\defu h{\theta}
\defu i{\infty}
\defu k{\kappa}
\defu l{\lambda}
\defu m{\mu}
\defu n{\nu}
\defu o{\omega}
\defu p{\pi}
\defu q{\psi}
\defu r{\rho}
\defu s{\sigma}
\defu t{\tau}
\defu w{\omega}
\defu x{\xi}
\defu y{\eta}
\defu z{\zeta}

\defu A{\forall}	
\defu C{{\noexpand\mathbb C}}
\defu D{\Delta}
\defu E{\exists}	
\defu F{\Phi}
\defu G{\Gamma}
\defu H{\Theta}
\defu L{\Lambda}
\defu O{\Omega}
\defu P{\Psi}
\defu Q{{\noexpand\mathbb Q}}
\defu R{{\noexpand\mathbb R}}
\defu S{\Sigma}
\defu Z{{\noexpand\mathbb Z}}

\ExtOn


\def\newterm#1{{\usefont{OT1}{cmss}{sbc}{n}#1}}

\def\FM{Fulton and MacPherson}
\def\Hopo{Hodge polynomial}
\def\Popo{Poincar\'e polynomial}

\def\B{{\mathrm B}}	
\def\C{{\mathbb C}}	
\def\F{{\mathcal F}}	
\def\I{{\mathcal I}}	
\def\L{{\mathcal L}}	
\def\O{{\mathcal O}}	
\def\P{{\mathbb P}}	
\def\S{{\mathcal S}}	
\def\x{{\mathbf x}}	

\def\Am{{\mathbb A\!^m}}	
\def\Bl{\operatorname{Bl}}	
\def\Conf{\mathrm F}		
\def\CS{{\mathrm C\mathrm S}}	
\def\Ln{{L_{[n]}}}		
\def\PS{{\mathrm P\mathrm S}}	
\def\Sym{\mathbb S}		
\def\kt{\Bbbk^\times}		

\def\0{{{}^\circ\!}}		
\def\1{{`p_1}}
\def\2{{`p_2}}
\def\<{\leqslant}
\def\>{\geqslant}
\def\suchthat{\thinspace|\thinspace}
\def\qed{\ \hfill$\square$}	

\def\pd<#1>{{\!\left<#1\right>}}	
\def\Xn{{X\pd<n>}}			
\def\normal#1#2{N_{#1\!/#2}}		
\def\mydot{\scriptstyle\bullet}		

\def\Bigfactor#1#2{
	\raisebox{0.5ex}{$#1$}\!\Big/\!\raisebox{-0.5ex}{$#2$}
}							
\def\BBigfactor#1#2{\Big(\Bigfactor{#1}{#2}\Big)}	

\def\mycaption#1#2{
	\addtocounter{figure}{1}
	\parbox[b]{#1}{{\scshape Figure $\arabic{figure}$.}\ #2}
}	
\def\mycaptionempty{
	\addtocounter{figure}{1}
	\parbox[b]{1in}{\scshape Figure $\arabic{figure}$}
}


\title[Polydiagonal compactification]{
Polydiagonal compactification
of~configuration spaces
}

\author[A. Ulyanov]{Alexander P. Ulyanov}
\address{
\noindent Department of Mathematics, The Pennsylvania State University
\newline 218 McAllister Building, University Park, PA 16802
}
\email{ulyanov@math.psu.edu}
\thanks{Research partially supported by NSF grant DMS--9803593}

\date{May 2, 2000}

\begin{abstract}
A smooth compactification $\Xn$ of the configuration space
of $n$ distinct labeled points in a smooth algebraic variety~$X$
is constructed by a~natural sequence of blowups,
with the full symmetry of the permutation group~$\Sym_n$
manifest at each stage.
The strata of the normal crossing divisor at infinity
are labeled by {\em leveled} trees and
their structure is studied.
This is the maximal wonderful compactification
in the sense of De Concini--Procesi,
and it has a strata-compatible surjection onto
the Fulton--MacPherson compactification.
The degenerate configurations added in the compactification
are geometrically described by {\em polyscreens}
similar to the screens of Fulton and MacPherson.

In characteristic~0, isotropy subgroups of
the action of~$\Sym_n$ on~$\Xn$ are abelian,
thus $\Xn$ may be a step toward
an explicit resolution of singularities of
the symmetric products~$X^n\!/\Sym_n$.
\end{abstract}

\maketitle

\section*{Introduction}

The configuration space $\Conf(X,n)$ of $n$ distinct labeled points
in a topological space~$X$ is
the complement in the Cartesian product~$X^n$
of the union of the large diagonals
$`D^{ij}=\{(x_1,\dots,x_n)\suchthat x_i=x_j\}$.
Pioneering studies of these spaces by
Fadell, Neuwirth, Arnold and Cohen 
\cite{Arnold,FCohen,Fadell,FadellNeuwirth}
evolved into a still active area of algebraic topology;
Totaro opens his paper with a brief review~\cite{Totaro}.
Somewhat later,
a compactification of $\Conf(\C,n)$ modulo affine automorphisms,
known as the Grothendieck--Knudsen moduli space of
stable $n$-pointed curves of genus~0,
rose to prominence in modern algebraic geometry
\cite{Deligne:resume,Kapranov:Veronese,Keel,Knudsen}.

Then \FM\ devised a powerful construction that
works for any nonsingular algebraic variety
and produces a compactification $X[n]$ of $\Conf(X,n)$
with a remarkable combination of properties~\cite{FM}:
\begin{list}{$\triangleright$}{}
\item
$X[n]$ is nonsingular.
\item
$X[n]$ naturally comes equipped with a proper map onto $X^n$.
\item
$X[n]$ is symmetric:
it carries an action of the symmetric group~$\Sym_n$
by permuting the labels.
\item
The complement $D=X[n]\smallsetminus \Conf(X,n)$ is
a normal crossing divisor.
\item
The combinatorial structure of~$D$
and of the resulting stratification of~$X[n]$
is explicitly described:
the components of~$D$ correspond to
the subsets of $[n]=\{1,\dots,n\}$ with at least~2 elements;
their intersections, the strata,
correspond to nested collections of such subsets,
and the latter are just a reincarnation of
rooted trees with~$n$ marked leaves.
\item
Degenerate configurations have simple geometric descriptions.
\end{list}

Further results of \FM\ include:
a functorial description of~$X[n]$,
used to prove many of its properties listed above;
a fact that all isotropy subgroups of~$\Sym_n$ acting on~$X[n]$
are solvable;
some intersection theory, namely,
a presentation of the intersection rings of $X[n]$ and of its strata,
and, as an application,
a computation of the rational cohomology ring of $\Conf(X,n)$
for~$X$ a smooth compact complex variety.

About the same time,
constructions related to the Fulton--MacPherson compactification appeared,
all motivated by, and suited to,
some problems of mathematical physics:
for real manifolds \cite{AxelrodSinger,Kontsevich:Feynman};
for complex curves~\cite{BeilinsonGinzburg},
with later extension to higher dimensions~\cite{Ginzburg}.

The compactifications $X[n]$ are defined inductively,
with the step from $X[n]$ to $X[n+1]$
performed by a sequence of blowups
$$
\xymatrix @+3pt
{
X[n+1]=
Y_{n}\ar[r]^<(0.3){{}\alpha_{n-1}} &
Y_{n-1}\ar[r]^<(0.35){{}\alpha_{n-2}} &
{}\dots\ar[r]^{{}\alpha_1} &
Y_1    \ar[r]^<(0.2){{}\alpha_0} &
Y_0=X[n]\times X,
}
$$
where the center of the blowup~$`a_k$
is a disjoint union of subvarieties in~$Y_k$ 
corresponding in a specified way to
the subsets of~$[n]$ of cardinality $n-k$.
Thus, the symmetry of $\Sym_{n+1}$ is not present
at the intermediate stages.
An alternative, and completely symmetric, description of~$X[n]$ as
the closure of $\Conf(X,n)$ in a product of blowups
does not provide much insight into the structure of~$X[n]$,
so the inductive sequence of blowups is essential for that.

Fulton and MacPherson remark:
\begin{quote}
{\small
It would be interesting to see
if other sequences of blowups give compactifications that are symmetric,
and whose points have explicit and concise descriptions
\cite[bottom of p.\! 196]{FM}.
}
\end{quote}

An example of such a compactification,
for any nonsingular algebraic variety~$X$,
is studied in the present paper.
I~denote it by~$\Xn$ and call it a {\em polydiagonal compactification},
because the blowup loci are not only the diagonals of~$X^n$,
but also their {\em intersections}.
The idea is very simple:
one who tries to blow up all diagonals of the same dimension simultaneously
is forced to blow up all their intersections prior to that,
and this prescribes the sequence.
Following \FM's terminology,
$\Xn$ is a compactification even though
it is only compact when~$X$ itself is compact.
In general, it is equipped with a canonical proper map onto~$X^n$.

This construction applies also to real manifolds,
with real blowups replacing algebraic blowups.
The compactification is then a manifold with corners,
and the results about the strata presented here
can be rephrased to describe the combinatorics of its boundary.

The construction of~$\Xn$ is in some respects
similar to that of~$X[n]$,
with one important difference:
the former is completely symmetric {\em at each stage}.
This reduces logical complexity of the construction
even though it involves (considerably) more blowups.
From this last fact stems another feature of~$\Xn$:
it distinguishes some collisions that are treated as equal by \FM.
There is a surjection $\vartheta_n \colon \Xn\to X[n]$ that
essentially retreats from making these distinctions,
and it is completely symmetric as well.
Regardless of~$X$,
this map,
derived from a description of~$\Xn$ as
the closure of $\Conf(X,n)$ in a product of blowups,
is an isomorphism for $n\<3$ only,
and an iterated blowup otherwise.
The fibers of~$\vartheta_n$ have purely com\-bi\-na\-to\-ri\-al nature
and do not depend even on the dimension of~$X$;
their detailed description will appear in a separate paper~\cite{U2}.

Geometrically, the limiting configurations in
the Fulton--MacPherson compactification are viewed in terms of
tree-like successions of {\em screens},
each of which is a tangent space to~$X$ with several labeled points in it,
considered modulo translations and dilations.
In a similar visualization for points in~$\Xn$,
labels of a new kind,
necessary because~$\Xn$ has \inquotes{more} points than~$X[n]$,
augment the screens.
This rests on a study of the strata:
they are bundles over $X\pd<r>$,
$r < n$,
with fibers decomposable into products of certain projective varieties.
Named {\em bricks},
they form a family indexed by integer partitions that includes,
for example,
permutahedral varieties.
The latter in fact show up in each brick as constituents
that account for those new labels.

As for the combinatorics underlying~$\Xn$,
here the place of subsets of~$[n]$,
nested collections of such subsets,
and plain rooted trees
is taken by partitions of the set~$[n]$,
chains of such partitions,
and rooted trees whose vertices are assigned integer numbers,
called {\em levels}.
With these changes,
the natural stratification of~$\Xn$ is
quite similar to that of~$X[n]$;
moreover,
$\vartheta_n$~is a strata-compatible map corresponding to
the forgetful map from leveled trees to usual rooted trees.

Analogues for~$\Xn$ of most results of \FM\ follow purely geometrically.
Since the proofs do not require a functorial description of the space,
it is omitted.

The action of the symmetric group~$\Sym_n$ on~$X^n$
by permuting the labels has fixed points.
Fulton and MacPherson showed that the isotropy subgroups
of the label permutation action of~$\Sym_n$ on~$X[n]$ are solvable
\cite[Theorem~5]{FM}.
It turns out that in characteristic~0
the similar action of~$\Sym_n$ on~$\Xn$ has
only abelian isotropy subgroups;
thus,
singularities of $\Xn\!/\Sym_n$
can in principle be resolved by toric methods 
\cite{Ash+,JLB:eventails,Kempf+,Oda}.
The resulting space will provide
an explicit desingularization of the symmetric product $X^n\!/\Sym_n$,
as well as a smooth compactification of $\B(X,n)=\Conf(X,n)/\Sym_n$,
the configuration space of~$n$ {\em unlabeled} points in~$X$.

De Concini and Procesi developed a general approach
to compactifying complements of linear subspace arrangements
by iterated blowups \cite{DeCP}.
For each arrangement,
it yields a family of {\em wonderful}\/ blowups
with minimal and maximal elements.
Although they work with linear subspaces,
their technique is local and can be applied to $X^n \smallsetminus \Conf(X,n)$
for any smooth variety~$X$;
in this case,
the Fulton--MacPherson compactification is the minimal one,
while the polydiagonal compactification is the maximal one.
Along the lines of De Concini, MacPherson and Procesi~\cite{MacPP},
Yi Hu has extended many results presented here in Sections
\ref{sec-construct},~\ref{sec-hodge} and~\ref{sec-closure}
to blowups of arrangements of smooth subvarieties and then recovered
Kirwan's partial desingularization
of geometric invariant theory quotients~\cite{Hu,Kirwan}.

In addition,
Hu computed the intersection rings in that general context of arrangements.
In the case of $\Xn$ these rings may be used to build
a differential graded algebra model of $\Conf(X,n)$
for~$X$ a compact complex algebraic manifold,
as \FM\ did.
After that,
Kriz streamlined their differential graded algebra,
while Totaro extracted a presentation
of the cohomology ring of the configuration space
from the Leray spectral sequence of its embedding 
into its \inquotes{naive}\ compactification~$X^n$~\cite{Kriz,Totaro}.

\subsection*{Historical note}
(Communicated by W.~Fulton.)
\FM\ sought to build the space
whose points would be described by screens;
early attempts led them to consider the spaces
denoted here by $X\pd<4>$ and~$X\pd<5>$,
and to identify what to blow down to create the desired $X[4]$ and~$X[5]$.
Seeing that as $n$ grows,
the blowdown description quickly becomes unwieldy,
they chose not to pursue this in general and
finally settled on a nonsymmetric procedure.
D.~Thurston pointed out a symmetric construction of $X[n]$
and used its real analogue in his work on knot invariants~\cite{Thurston}.

\subsection*{Standing assumptions}
Throughout the paper,
$X$ is a smooth irreducible $m$-dimensional ($m > 0$)
algebraic variety over some field~$\Bbbk$,
and $n$ is the number of labeled points in~$X$.
The section on Hodge polynomials applies only to complex varieties,
and that on the symmetric group action,
only to the characteristic~0 case.

\subsection*{Outline of the paper}
The first section is informal
and serves to introduce the basic ideas of the polydiagonal compactification
on the simplest example.
A combinatorial interlude of Section~\ref{sec-combinat}
is followed by a discussion of polyscreens and colored screens
that represent points in~$\Xn$.

Formally stated and proved results begin
in Section~\ref{sec-construct} that contains:
construction of $\Xn$ by a symmetric sequence of blowups,
a description of the combinatorics of
the complement $\Xn\smallsetminus\Conf(X,n)$
as a divisor with normal crossings
and of the ensuing stratification of~$\Xn$,
and a recurrent formula for the number of the strata.
If~$X$ is a complex variety,
the blowup construction translates into a formula for
the (virtual) \Hopo\ $e(\Xn)$ in terms of~$e(X)$
derived in the next section.
In Section~\ref{sec-closure},
a consideration of~$\Xn$ as the closure of $\Conf(X,n)$ in a product of blowups
implies a surjection $\vartheta_n \colon \Xn\to X[n]$,
written then as an iterated blowup.
Technical analysis of the strata of~$\Xn$ occupies Section~\ref{sec-strata},
and the last section deals with
the isotropy subgroups of~$\Sym_n$ acting on~$\Xn$.

\subsection*{Acknowledgements}
In many ways,
I~am indebted to Jean--Luc Brylinski,
my Ph.D.~advisor.
I~am most grateful to William Fulton and Jim Stasheff
for sharing their advice and for many valuable comments and insights.
The referee's suggestions helped me improve exposition.
I~would also like to thank Dmitry Tamarkin for useful discussions.

\section{Small numbers of colliding points}
\label{sec-x4}

The purpose of this section is
to introduce the main ideas of the paper
by looking at the case of 4 points---%
the smallest integer~$n$ for which $\Xn$ is different from $X[n]$ is~4.

To begin with,
consider an example of two collisions of four points in $X=`C^2$.
The corresponding two limiting configurations
arising in the approach of \FM\ coincide;
however, 
the polydiagonal compactification will distinguish them.
Take four points labeled by 1 through 4
and make them collide as $t\to 0$ in the following way:
\begin{list}{$\diamond$}{}
	\item the distance between 1 and 2 is $O(t^3)$,
	\item the distance between 3 and 4 is $O(t^2)$,
	\item the distance between the two pairs (12) and (34) is $O(t)$.
\end{list}
Then do the same thing, except for a small exchange:
\begin{list}{$\diamond$}{}
	\item the distance between 1 and 2 is $O(t^2)$,
	\item the distance between 3 and 4 is $O(t^3)$,
\end{list}
and call the two limiting points $\x_1$ and $\x_2$.

Both limiting points lie in the same stratum of~$X[4]$,
the intersection of three divisors $D(1234)$, $D(12)$, and $D(34)$.
The dimension of this stratum is~5;
the dimension of its fiber over a point in the small diagonal
$`D\subset X^4$ is~3.
The three parameters record
the \inquotes{directions}\ 
of collisions encoded by the middle tree in Figure~\ref{fig-not-same}.
Specifying these directions for the two approach curves,
that is,
vectors hidden behind the symbol~$O$,
one can arrange that $\x_1=\x_2$ in $X[4]$.

\begin{figure}[b]
\label{fig-not-same}
$$
\xy
(-30,0)*{\ };
(94,2)*{\parbox[b]{1.6in}{\caption{}}};
(-6,-6.4)*{{}^1},
(-2,-6.4)*{{}^2},
(2,-6.4) *{{}^3},
(6,-6.4) *{{}^4},
@={(-6,-4),(-4,0),(-2,-4),(-4,0),(0,8),(2,4),(2,-4),(2,4),(6,-4)},
s4="prev" @@{;"prev";**@{-}="prev"};
(24,-6.4)*{{}^1},
(28,-6.4)*{{}^2},
(32,-6.4) *{{}^3},
(36,-6.4) *{{}^4},
@i; @={(24,-4),(27,2),(28,-4),(27,2),(30,8),(33,2),(32,-4),(33,2),(36,-4)},
s4="prev" @@{;"prev";**@{-}="prev"};
(54,-6.4)*{{}^1},
(58,-6.4)*{{}^2},
(62,-6.4) *{{}^3},
(66,-6.4) *{{}^4},
@i; @={(54,-4),(58,4),(58,-4),(58,4),(60,8),(64,0),(62,-4),(64,0),(66,-4)},
s4="prev" @@{;"prev";**@{-}="prev"};
(22,2) \ar @{~>} (38,2);
(52,2) \ar @{~>} (8,2)
\endxy
$$
\end{figure}

These approach curves actually belong to
a whole family~$\mathcal F$ of curves in $\Conf(X,4)$
whose limits in $X[4]$ may coincide.
Indeed,
consider the diagonals $`D^{12}$ and $`D^{34}$ in $X^4$,
and their intersection $`D^{12|34}$.
Both curves approach this intersection,
but the first one does it
while having a 3rd degree osculation to $`D^{12}$,
and the second one does the same with $`D^{34}$.
The projectivized normal space $\P(T_p X^4/T_p `D^{12|34})$
parametrizes the family,
and the two curves above correspond to
normal directions going along $`D^{12}$ and $`D^{34}$ respectively.
This suggests looking into a possibility of involving
blowups of subvarieties like $`D^{12} \cap `D^{34}$,
if the objective is to obtain a compactification that
would distinguish from one another collisions
produced by curves in such families.

The space that achieves this
results from implementing
a simple idea of blowing up \inquotes{from the bottom to the top}.
Although the dominant feature of the general case
first comes to light when $n=4$,
it may be useful to begin with the cases of two and three colliding points.

Assume that $\dim X > 1$.
There is no ambiguity about the case of $n=2$ points:
the compactification is the blowup of the diagonal in~$X^2$.
If~$n=3$, blowing up the small diagonal $`D \subset X^3$
creates disjoint proper transforms of $`D^{12}$, $`D^{13}$ and~$`D^{23}$
that can then be blown up in any order.
The resulting compactification coincides with~$X[3]$.
For~$n>3$, however, this strategy will not work,
and some additional blowups are needed \cite[bottom of p.~196]{FM},
but what they are \FM\ do not specify.

The left graph in Figure~3 
shows the diagonals in~$X^4$,
including the space itself,
as vertices,
and (nonrefinable) inclusions of the diagonals into each other as edges.
As before,
blow up the small diagonal first,
then blow up the (disjoint) proper transforms of
the four larger diagonals,
like $`D^{123}$.
Now try to blow up the next level below them simultaneously.
It does not work:
these six largest diagonals have not been made disjoint.
How can this be fixed?

\begin{figure}[t]
\label{fig-six-lines}
\def\0{\kern-1pt}
\xy
(-40,10)*{\ };
(0,-8)*{\ };
(46,10)*{\parbox[b]{1.6in}{\caption{}}};
(0,-6); (0,18)**\dir{-};
(-14,0); (14,0)**\dir{-};
(11.4,-5.1); (-3.6,17.4)**\dir{-};
(-11.4,-5.1); (3.6,17.4)**\dir{-};
(14,-3); (-10,9)**\dir{-};
(-14,-3); (10,9)**\dir{-};
(0,0)*{\mydot};
(0,4)*{\mydot};
(8,0)*{\mydot};
(4,6)*{\mydot};
(0,12)*{\mydot};
(-4,6)*{\mydot};
(-8,0)*{\mydot};
(3.9,15.1)*{\scriptstyle 1\!2};
(9,-4)*{\scriptstyle 1\!3};
(-12.6,1.2)*{\scriptstyle 1\!4};
(-1.4,-4.2)*{\scriptstyle 2\03};
(9.5,7.2)*{\scriptstyle 2\04};
(-8.2,9.7)*{\scriptstyle 3\04};
\endxy
\end{figure}

\begin{figure}[b]
\label{fig-x4}
\def\0{\kern-1pt}
\hskip3mm
\xy
(-28,14)*{\ };
(60,-18)*{
	\mycaption{6in}{Diagonals (left) and polydiagonals (right) in~$X^4$}};
%
%
(0,12)="1234",
(-18,4)="123",
(-6,4)="124",
(6,4)="134",
(18,4)="234",
(-15,-4)="12",
(-9,-4)="13",
(-3,-4)="23",
(3,-4)="14",
(9,-4)="24",
(15,-4)="34",
(0,-12)="no",
@={"1234","123","12","124","24","234","34","134","1234",
"234","23","123","13","134","14","124","1234"},
s0="prev" @@{;"prev";**@{-}="prev"},
@i @={"12","13","23","14","24","34"},
"no"; @@{**@{-}},
@i @={"1234","123","124","134","234","12","13","14","23","24","34","no"},
@@{*{\mydot}},
"1234"+(0,1)*{{}^{1\02\03\04}},
"123"+(0,1)*{{}^{1\02\03}},
"124"+(-.5,1)*{{}^{1\02\04}},
"134"+(.5,1)*{{}^{1\03\04}},
"234"+(0,1)*{{}^{2\03\04}},
"12"+(-2,-1)*{{}^{1\02}},
"13"+(-2,-1)*{{}^{1\03}},
"23"+(-2,-1)*{{}^{2\03}},
"14"+(2,-1)*{{}^{1\04}},
"24"+(2,-1)*{{}^{2\04}},
"34"+(2,-1)*{{}^{3\04}},
%
%
(70,0)="m",
"m"+(0,12)="1234",
"m"+(-6,5)="12-34",
"m"+(0,5)="13-24",
"m"+(6,5)="14-23",
"m"+(-18,4)="123",
"m"+(-12,4)="124",
"m"+(12,4)="134",
"m"+(18,4)="234",
"m"+(-15,-4)="12",
"m"+(-9,-4)="13",
"m"+(-3,-4)="23",
"m"+(3,-4)="14",
"m"+(9,-4)="24",
"m"+(15,-4)="34",
"m"+(0,-12)="no",
@i @={"1234","123","12","124","24","234","34","134","1234",
"234","23","123","13","134","14","124","1234"},
s0="prev" @@{;"prev";**@{-}="prev"},
@i @={"12-34","13-24","14-23"},
"1234"; @@{**@{-}},
@i @={"no","12","12-34","34","no","23","14-23","14","no",
"13","13-24","24","no"},
s0="prev" @@{;"prev";**@{-}="prev"},
@i @={"1234","123","124","134","234","12-34","13-24","14-23",
"12","13","14","23","24","34","no"},
@@{*{\mydot}},
"1234"+(0,1)*{{}^{1\02\03\04}},
"123"+(0,1)*{{}^{1\02\03}},
"124"+(-.5,1)*{{}^{1\02\04}},
"134"+(.5,1)*{{}^{1\03\04}},
"234"+(0,1)*{{}^{2\03\04}},
"12"+(-2,-1)*{{}^{1\02}},
"13"+(-2,-1)*{{}^{1\03}},
"23"+(-2,-1)*{{}^{2\03}},
"14"+(2,-1)*{{}^{1\04}},
"24"+(2,-1)*{{}^{2\04}},
"34"+(2,-1)*{{}^{3\04}},
\endxy
\end{figure}


The six lines intersecting at seven points
depicted in Figure~2 
are the images of the large diagonals of~$`R^4$
in the real projective plane $\P(`R^4/`D)$,
where~$`D$ is the small diagonal.
Four of the points correspond to diagonals like $`D^{123}$,
and the other three,
where the intersections are normal,
represent additional loci that need to be blown up
to make the large diagonals disjoint.
The second graph in Figure~3 is obtained from the first one
by adding these three intersections
$`D^{12} \cap `D^{34}$, $`D^{13} \cap `D^{24}$ and $`D^{14} \cap `D^{23}$.
All seven vertices in the second row correspond to subvarieties pairwise
disjoint after the blowup of the small diagonal $`D\subset X^4$,
so they can be blown up simultaneously,
and---crucially---after that the subvarieties from the row just below
become disjoint and can be blown up simultaneously.
This gives a compactification $X\pd<4>$ of~$\Conf(X,4)$.

\begin{figure}[t]
\label{fig-FM-12-34}
\xy
	(30,-7)*{\mycaption{3in}{A point in $X[4]$}};
	(-5,29); (18,29)**\crv{(-15,18)&(4,20)&(10,26)&(30,18)&(26,26)};
	(-2,25)*{X};
	(12,26)*{\mydot};
	@={(8,12),(8,22),(18,22),(18,12)};
	s0="prev" @@{;"prev";**@{-}="prev"};
	(8,22); (11.5,26.1)**\dir{.};
	(18,22); (12.3,26.2)**\dir{.};
	(14.5,15)*{\mydot};
	(11.5,19)*{\mydot};
    @i;	@={(0,0),(0,10),(10,10),(10,0)};
	s0="prev" @@{;"prev";**@{-}="prev"};
	(0,10); (11.2,19.4)**\dir{.};
	(10,0); (11.9,19)**\dir{.};
	(2,6)*{\mydot} -(0,1.3)*{\scriptscriptstyle 1};
	(8,4)*{\mydot} -(0,1.3)*{\scriptscriptstyle 2};
    @i;	@={(16,0),(16,10),(26,10),(26,0)};
	s0="prev" @@{;"prev";**@{-}="prev"};
	(16,0); (14.1,15)**\dir{.};
	(26,10); (14.7,15.4)**\dir{.};
	(18,3)*{\mydot} -(0,1.3)*{\scriptscriptstyle 3};
	(24,7)*{\mydot} -(0,1.3)*{\scriptscriptstyle 4};
	(70,2)="m"
	+(33,-8.5)*{\mycaption{2in}{Dilation}};
    @i; @={"m"+(0,0),"m"+(16,0),"m"+(16,16),"m"+(0,16)};
	s0="prev" @@{;"prev";**@{-}="prev"};
	"m"+(5,6)*{\mydot}	+(1,-1)*{\scriptscriptstyle 3};
	"m"+(11,10)*{\mydot}	-(1,-1)*{\scriptscriptstyle 4};
	"m"+(11.7,-1.8)*{\scriptstyle\alpha_{34}=2};
	"m"+(30,0)="m";
    @i; @={"m"+(0,0),"m"+(16,0),"m"+(16,16),"m"+(0,16)};
	s0="prev" @@{;"prev";**@{-}="prev"};
	"m"+(2,4)*{\mydot}	+(1,-1)*{\scriptscriptstyle 3};
	"m"+(14,12)*{\mydot}	-(1,-1)*{\scriptscriptstyle 4};
	"m"+(12,-1.8)*{\scriptstyle\alpha_{34}=1};
	"m"+(-7,8)*{=};
\endxy
\end{figure}

The construction of $X\pd<4>$ involves three more blowups than that of $X[4]$,
so the complement of $\Conf(X,4)$ in $X\pd<4>$ has three additional components
$D^{12|34}$, $D^{13|24}$ and $D^{14|23}$.
Collisions belonging to the family~$\mathcal F$ discussed above
result in points in $Z = D^{1234}\cap D^{12|34}$.
To accommodate these,
as well as more complicated degenerations of the same nature
that appear for $n>4$,
two new features are added to Fulton--MacPherson screens:
the screens are grouped into \newterm{levels},
and the group on each level bears
a new parameter living in a projective space.

Figure~4 illustrates the screen description
of the limiting points in $X[4]$ of the family~$\mathcal F$:
its macroscopic part is a single point in~$X$ and
its microscopic part consists of three screens,
one for each of the subsets $1234$, $12$~and~$34$ of $\{1,2,3,4\}$.
A~screen is a tangent space $T_p X$
with a configuration of points in it,
considered up to dilations and translations.
In particular,
the last two screens,
$S_{12}$ and $S_{34}$,
are completely independent of each other.

\begin{figure}[b]
\label{fig-new-level}
\xy
	(30,-12)*{\mycaption{3in}{A degeneration in $X\pd<4>$}};
	(-5,29); (18,29)**\crv{(-15,18)&(4,20)&(10,26)&(30,18)&(26,26)};
	(-2,25)*{X};
	(12,26)*{\mydot};
	@={(8,12),(8,22),(18,22),(18,12)};
	s0="prev" @@{;"prev";**@{-}="prev"};
	(8,22); (11.5,26.1)**\dir{.};
	(18,22); (12.3,26.2)**\dir{.};
	(14.5,15)*{\mydot};
	(11.5,19)*{\mydot};
    @i;	@={(0,0),(0,10),(10,10),(10,0)};
	s0="prev" @@{;"prev";**@{-}="prev"};
	(0,10); (11.2,19.4)**\dir{.};
	(10,0); (11.9,19)**\dir{.};
	(2,6)*{\mydot} -(0,1.3)*{\scriptscriptstyle 1};
	(8,4)*{\mydot} -(0,1.3)*{\scriptscriptstyle 2};
	(6,-1.8)*{\scriptstyle\alpha_{12}=5};
    @i;	@={(16,0),(16,10),(26,10),(26,0)};
	s0="prev" @@{;"prev";**@{-}="prev"};
	(16,0); (14.1,15)**\dir{.};
	(26,10); (14.7,15.4)**\dir{.};
	(18,3)*{\mydot} -(0,1.3)*{\scriptscriptstyle 3};
	(24,7)*{\mydot} -(0,1.3)*{\scriptscriptstyle 4};
	(22,-1.8)*{\scriptstyle\alpha_{34}=1};
%
	(65,29); (88,29)**\crv{(55,18)&(74,20)&(80,26)&(100,18)&(96,26)};
	(68,25)*{X};
	(82,26)*{\mydot};
    @i;	@={(78,12),(78,22),(88,22),(88,12)};
	s0="prev" @@{;"prev";**@{-}="prev"};
	(78,22); (81.5,26.1)**\dir{.};
	(88,22); (82.3,26.2)**\dir{.};
	(84.5,15)*{\mydot};
	(81.5,19)*{\mydot};
    @i;	@={(70,0),(70,10),(80,10),(80,0)};
	s0="prev" @@{;"prev";**@{-}="prev"};
	(70,10); (81.2,19.4)**\dir{.};
	(80,0); (81.9,19)**\dir{.};
	(72,6)*{\mydot} -(0,1.3)*{\scriptscriptstyle 1};
	(78,4)*{\mydot} -(0,1.3)*{\scriptscriptstyle 2};
    @i;	@={(86,0),(86,10),(96,10),(96,0)};
	s0="prev" @@{;"prev";**@{-}="prev"};
	(86,0); (84.1,15)**\dir{.};
	(96,10); (84.7,15.4)**\dir{.};
	(91,5)*{\mydot}
    @i;	@={(86,-12),(86,-2),(96,-2),(96,-12)};
	s0="prev" @@{;"prev";**@{-}="prev"};
	(86,-2); (90.7,5.2)**\dir{.};
	(96,-2); (91.3,5.2)**\dir{.};
	(88,-9)*{\mydot} -(0,1.3)*{\scriptscriptstyle 3};
	(94,-5)*{\mydot} -(0,1.3)*{\scriptscriptstyle 4};
	(110,25)*{\txt{Level 0}};
	(110,15)*{\txt{Level 1}};
	(110, 5)*{\txt{Level 2}};
	(110,-5)*{\txt{Level 3}};
	(47,5)*\txt{as $\alpha_{34} \to 0$};
        (60,7) \ar @/^0.4ex/ (35,7);
\endxy
\end{figure}

Pictures like the left one in Figure~6 will represent generic points of~$Z$.
It~consists of three levels:
\begin{enumerate}
\renewcommand\theenumi{\arabic{enumi}}
\setcounter{enumi}{-1}
\item	one point in~$X$,
\item	a screen for $1234$ with two distinct points,
	and
\item	a pair of screens $S_{12}$ and $S_{34}$
	together with their \newterm{scale factors} $`a_{12}$ and $`a_{34}$,
	where the pair $[`a_{12} \!:\! `a_{34}]$
	is considered as a point in~$\P^1$.
\end{enumerate}
The scale factors serve to compare the approach speeds of the pairs 12~and~34
by keeping track of independent dilations of their respective screens:
for all nonzero scalars~$`f_i$,
the pairs $(S_i,`a_i)$ and $(`f_i S_i,`a_i/`f_i)$ are identified,
where the screen in the second pair
is the dilation of~$S_i$ by the factor of~$`f_i$,
as in Figure~5.


Nongeneric points of~$Z$,
which lie in $Z\cap D^{12}$ and $Z\cap D^{34}$,
correspond to incomparable speeds
and to the points $[0\!:\!1]$ and $[1\!:\!0]$ in~$\P^1$.
They result from collisions mentioned
in the beginning of the section.
Keeping the screen~$S_{34}$ fixed while letting $`a_{34} \to 0$
is the same thing as
keeping $`a_{34}$ fixed while contracting the screen.
In the limit the two points in it collide,
but a new screen appearing on level~$3$ separates them.
Trivial screens,
which contain a single point,
may be omitted from the pictures.

Similarly,
points in $D^{12|34}$ away from $D^{1234}$
are represented by configurations of two distinct points in~$X$,
labeled 12 and 34,
plus screens $S_{12}$ and $S_{34}$ together with their scale factors,
generically on the same level and degenerating to two levels.

The microscopic levels in Figure~6 correspond to the intersecting divisors:
the first to $D^{1234}$,
the second to $D^{12|34}$
and,
in the right half of the figure,
the third to $D^{34}$.
Accordingly,
trees that link screens together acquire some extra structure:
levels of vertices.
For example,
the two pictures in Figure~6 correspond to
the middle and right trees in Figure~1.
Such trees index the strata of~$X\pd<4>$.

Scale factors are redundant on any level that
contains only one nontrivial screen.
Since the middle tree in Figure~1 is,
up to relabeling,
the only tree with four leaves
in which two vertices may be on the same level,
points in $X\pd<4>$ outside the three additional divisors
will have exactly the same screen description as for~$X[4]$.
In fact,
forgetting the scale factors gives
a map~$\vartheta_4 \colon X\pd<4> \to X[4]$
that blows down the divisor $D^{12|34}$ to the stratum $D(12) \cap D(34)$,
and respectively for $D^{13|24}$ and $D^{14|23}$.

A combinatorial basis is necessary
in order to generalize these ideas to an arbitrary number of points,
and it is very easy to find.
The definition of~$`D^{S}$ for any subset~$S$ of $[n]=\{1,\dots,n\}$
applied to $S=\{k\}$ gives $`D^{\{k\}}=X^n$,
hence $`D^{123}=`D^{123} \cap `D^{4}$ and so on.
The true combinatorial basis will thus be
{\em the partitions of the set\/}~$[n]$.
Indeed,
when $n=4$,
the first blowup is that of $`D=`D^{1234}$,
which corresponds to the only partition into one block;
the next stage blowup centers correspond to all partitions into two blocks;
finally, all those corresponding to partitions into three blocks
are blown up:
$`D^{12}=`D^{12} \cap `D^{3} \cap `D^{4}$ and so on.

\section{Combinatorial background}
\label{sec-combinat}

This section is a short primer on the language of the rest of the paper:
it deals with basic properties of set partitions
and a bijection between partition chains and leveled trees.


Let $[n]$ denote the set $\{1,\dots,n\}$ of integers.
A~partition~$`p$ of~$[n]$ is a set of disjoint subsets of~$[n]$,
called the blocks of~$`p$,
whose union is~$[n]$.
Nonsingleton blocks are called~\newterm{essential}.
The two functions of partitions that are most important for this work
are $`r(`p)$,
the number of blocks,
and $`e(`p)$,
the number of essential blocks.
The integer partition whose parts are {\em one less than}
the cardinalities of the essential blocks of~$`p$
is called the \newterm{essential shape} of~$`p$ and denoted by~$`l(`p)$.
For example,
$\1 = \{12357,9,468\}$ and $\2 = \{15,23,7,9,468\}$
are two partitions of~$[9]$
with
\begin{equation*}
\begin{align*}
{}`r(\1) &= 3,		& {}`e(\1) &= 2,	& {}`l(\1) &= (4,2),	&&\\
{}`r(\2) &= 5,		& {}`e(\2) &= 3, 	& {}`l(\2) &= (2,1,1).	&&
\end{align*}
\end{equation*}

Let~$\Ln$ be the set of all partitions of $[n]$.
There is a refinement partial order on~$\Ln$:
$\1\leqslant\2$ whenever each block of~$\2$ is contained in a block of~$\1$,
as in the example.
This makes $\Ln$ a ranked lattice,
with $`r(`p)$ being the rank function.
The minimal (bottom) and maximal (top) elements
of~$\Ln$ are denoted by $\bot$ and $\top$ respectively.

The Stirling number of the second kind $S(n,k)$
is the number of partitions of~$[n]$ into exactly~$k$ blocks.
Many textbooks on combinatorics discuss
these numbers and the partition lattice,
for instance,
Andrews~\cite{Andrews} and Stanley~\cite{Stanley}.

An interval $[`p',`p'']$ in a lattice~$L$ is
its subset $\{`p \suchthat `p' \< `p \< `p''\}$.
In~$\Ln$,
every lower interval $[\bot,`p]$ is isomorphic to $L_{[`r(`p)]}$ and
every upper interval $[`p,\top]$ is isomorphic to
$L_{[`n_1+1]} \times\dots\times L_{[`n_r+1]}$,
where $`l = `l(`p) = (`n_1,\dots,`n_r)$ is the essential shape of~$`p$.
This product will be denoted by~$L_`l$.


A totally ordered subset of a partially ordered set is called a chain.
The length of a chain is the number of its elements.
Half of the chains in $\Ln$ contain the top (finest) partition,
and the other half do not;
from now on,
a chain will mean a partition chain of the latter kind.

\begin{figure}[b]
\label{fig-leveled}
\xy
	(88,-8)*{\mycaption{7in}{
			From a partition chain to a leveled tree (and back)}};
	(-7,14.7)*{\scriptstyle \pi_1};
	(-7, 9.7)*{\scriptstyle \pi_2};
	(-7, 4.7)*{\scriptstyle \pi_3};
	(16,20)="root"		*{\bullet};
	(10,15)="a12357"	*{\scriptstyle 12357};
	(18,15)="a9"		*{\scriptstyle 9};
	(24,15)="a468"		*{\scriptstyle 468};
	(4,10)="b15"		*{\scriptstyle 15};
	(10,10)="b23"		*{\scriptstyle 23};
	(14.75,10)="b7"		*{\scriptstyle 7};
	(19.8,10)="b9"		*{\scriptstyle 9};
	(26.5,10)="b468"	*{\scriptstyle 468};
	(0.6,5)="c1"		*{\scriptstyle 1};
	(4.3,5)="c5"		*{\scriptstyle 5};
	(10,5)="c23"		*{\scriptstyle 23};
	(15.5,5)="c7"		*{\scriptstyle 7};
	(20,5)="c9"		*{\scriptstyle 9};
	(25.5,5)="c46"		*{\scriptstyle 46};
	(29.8,5)="c8"		*{\scriptstyle 8};
	(0,0)="d1"		*{\scriptstyle 1};
	(4,0)="d5"		*{\scriptstyle 5};
	(8,0)="d2"		*{\scriptstyle 2};
	(12,0)="d3"		*{\scriptstyle 3};
	(16,0)="d7"		*{\scriptstyle 7};
	(20,0)="d9"		*{\scriptstyle 9};
	(24,0)="d4"		*{\scriptstyle 4};
	(28,0)="d6"		*{\scriptstyle 6};
	(32,0)="d8"		*{\scriptstyle 8};
	"root";			"a12357"+(1,1.2)**\dir{-};
	"root";			"a9"+(-0.3,1.2)**\dir{-};
	"root";			"a468"+(-1.3,1.2)**\dir{-};
	"a12357"-(2,1.2);	"b15"+(1,1.2)**\dir{-};
	"a12357"-(0,1.2);	"b23"+(0,1.2)**\dir{-};
	"a12357"-(-2,1.2);	"b7"+(-0.5,1.2)**\dir{-};
	"a9"-(-0.5,1.2);	"b9"+(-0.3,1.2)**\dir{-};
	"a468"-(-0.5,1.2);	"b468"+(-0.3,1.2)**\dir{-};
	"b15"-(1,1.2);		"c1"+(0.3,1.2)**\dir{-};
	"b15"-(-0.5,1.2);	"c5"+(0,1.2)**\dir{-};
	"b23"-(0,1.2);		"c23"+(0,1.2)**\dir{-};
	"b7"-(-0.3,1.2);	"c7"+(0,1.2)**\dir{-};
	"b9"-(0,1.2);		"c9"+(0,1.2)**\dir{-};
	"b468"-(0.2,1.2);	"c46"+(0.2,1.2)**\dir{-};
	"b468"-(-1,1.2);	"c8"+(-0.5,1.2)**\dir{-};
	"c1"-(0.2,1.2);		"d1"+(0,1.2)**\dir{-};
	"c5"-(0.1,1.2);		"d5"+(0,1.2)**\dir{-};
	"c23"-(0.7,1.2);	"d2"+(0,1.2)**\dir{-};
	"c23"-(-0.7,1.2);	"d3"+(0,1.2)**\dir{-};
	"c7"-(-0.2,1.2);	"d7"+(0,1.2)**\dir{-};
	"c9"-(0,1.2);		"d9"+(0,1.2)**\dir{-};
	"c46"-(0.5,1.2);	"d4"+(0,1.2)**\dir{-};
	"c46"-(-1,1.2);		"d6"+(-0.3,1.2)**\dir{-};
	"c8"-(-0.5,1.2);	"d8"+(-0.2,1.2)**\dir{-};
(35,13)*{\leftrightarrow};
	(39,0)="m";
	"m"+(16,20)="root"	*{\bullet};
	"m"+(10,15)="12357"	*{\scriptstyle 12357};
	"m"+(26.5,10)="468"	*{\scriptstyle 468};
	"m"+(4,10)="15"		*{\scriptstyle 15};
	"m"+(10,5)="23"		*{\scriptstyle 23};
	"m"+(25.5,5)="46"	*{\scriptstyle 46};
	"m"+(0,0)="1"		*{\scriptstyle 1};
	"m"+(4,0)="5"		*{\scriptstyle 5};
	"m"+(8,0)="2"		*{\scriptstyle 2};
	"m"+(12,0)="3"		*{\scriptstyle 3};
	"m"+(16,0)="7"		*{\scriptstyle 7};
	"m"+(20,0)="9"		*{\scriptstyle 9};
	"m"+(24,0)="4"		*{\scriptstyle 4};
	"m"+(28,0)="6"		*{\scriptstyle 6};
	"m"+(32,0)="8"		*{\scriptstyle 8};
	"root";			"12357"+(1,1.2)**\dir{-};
	"root";			"468"+(-1,1.2)**\dir{-};
	"root";			"9"+(0,1.2)**\dir{-};
	"12357"-(1.5,1.2);	"15"+(1,1.2)**\dir{-};
	"12357"-(0,1.2);	"23"+(0,1.2)**\dir{-};
	"12357"-(-1.5,1.2);	"7"+(0,1.2)**\dir{-};
	"468"-(0.5,1.2);	"46"+(0,1.2)**\dir{-};
	"468"-(-0.5,1.2);	"8"+(-0.2,1.2)**\dir{-};
	"15"-(0.5,1.2);		"1"+(0,1.2)**\dir{-};
	"15"-(-0.5,1.2);	"5"+(0,1.2)**\dir{-};
	"23"-(0.7,1.2);		"2"+(0,1.2)**\dir{-};
	"23"-(-0.7,1.2);	"3"+(0,1.2)**\dir{-};
	"46"-(0.5,1.2);		"4"+(0,1.2)**\dir{-};
	"46"-(-1,1.2);		"6"+(-0.3,1.2)**\dir{-};
"m"+(35,13)*{\leftrightarrow};
	"m"+"root"-(0,0.4)="nroot";
	"m"+"12357"	="n12357";
	"m"+"468"	="n468";
	"m"+"15"+(0,0.4)="n15";
	"m"+"23"+(0,0.8)="n23";
	"m"+"46"+(0,0.8)="n46";
	"m"+"1"+(0,1.2)	="n1";
	"m"+"2"+(0,1.2)	="n2";
	"m"+"3"+(0,1.2)	="n3";
	"m"+"4"+(0,1.2)	="n4";
	"m"+"5"+(0,1.2)	="n5";
	"m"+"6"+(0,1.2)	="n6";
	"m"+"7"+(0,1.2)	="n7";
	"m"+"8"+(0,1.2)	="n8";
	"m"+"9"+(0,1.2)	="n9";
      @={"nroot","n12357","n15","n1","n15","n5","n15","n12357","n23","n2",
		"n23","n3","n23","n12357","n7","n12357","nroot","n9","nroot",
		"n468","n46","n4","n46","n6","n46","n468","n8","n468","nroot"};
s0="prev" @@{;"prev";**@{-}="prev"};
	"nroot"	*{\bullet};
	"n12357"*{\scriptscriptstyle\bullet};
	"n468"*{\scriptscriptstyle\bullet};
	"n15"*{\scriptscriptstyle\bullet};
	"n23"*{\scriptscriptstyle\bullet};
	"n46"*{\scriptscriptstyle\bullet};
	"n1"-(0,1.2)	*{\scriptstyle 1};
	"n2"-(0,1.2)	*{\scriptstyle 2};
	"n3"-(0,1.2)	*{\scriptstyle 3};
	"n4"-(0,1.2)	*{\scriptstyle 4};
	"n5"-(0,1.2)	*{\scriptstyle 5};
	"n6"-(0,1.2)	*{\scriptstyle 6};
	"n7"-(0,1.2)	*{\scriptstyle 7};
	"n8"-(0,1.2)	*{\scriptstyle 8};
	"n9"-(0,1.2)	*{\scriptstyle 9};
	(114,20)*\txt{\tiny level 0};
	(114,15.33)*\txt{\tiny level 1};
	(114,10.66)*\txt{\tiny level 2};
	(114, 6)*\txt{\tiny level 3};
\endxy
\end{figure}

Lengyel represented~\cite{Lengyel} partition chains as trees.
If $`g = \{`p_1,\dots,`p_k\}$,
where $`p_i < `p_{i+1}$ for $1 \< i \< k$,
then the associated tree has the blocks of each partition
as its interior vertices,
one additional vertex (the root)
and leaves labeled by $1,\dots,n$.
Edges indicate inclusions of blocks of $`p_{i+1}$ into those of~$`p_i$
and of the elements of~$[n]$ into the blocks of~$`p_k$;
they also connect the blocks of~$\1$ to the root.
The left tree in Figure~7 goes with the chain $`g = \{`p_1,`p_2,`p_3\}$,
where
$$
`p_1 = \{12357,9,468\}, \quad
`p_2 = \{15,23,7,9,468\}, \quad
`p_3 = \{1,5,23,7,9,46,8\}.
$$

The $2$-valent vertices
(except for the root if it happens to be such)
may be called the \newterm{phantom vertices}
because it is often convenient to omit them;
this gives trees like the middle one in the same figure.
Furthermore,
labels of interior vertices are also unnecessary.
In thus simplified tree the set of interior vertices
is the set $\{12357,468,23,46,15\}$
of all essential blocks in the three partitions,
and they appear to be on different {\em levels}
reflecting how far in the chain they survive unsubdivided.
This leads to the following

\begin{defn}
A \newterm{$k$-leveled tree} is a pair $(T,`y)$,
where~$T$ is a rooted tree without $2$-valent vertices,
except possibly for the root,
and~$`y$ is a surjective poset map from 
the set of vertices of~$T$ with the parent-descendant partial order
to the set of integers $\{0,\dots,k\}$ with its standard order.
(The root goes to~$0$.)
The number~$`y(v)$ is called the \newterm{level} of the vertex~$v$.
The map from leveled trees with marked leaves to
usual rooted trees with marked leaves by $(T,`y)\mapsto T$
is denoted by~$`h$.
\end{defn}

The term {\em leveled tree} belongs to Loday,
although his trees are binary~\cite{Loday}.
An~inspiring picture evinces that
Tonks used leveled trees implicitly~\cite{Tonks}.
In~both references the leaves are not marked.
The sole purpose of the root is to simplify wording:
without it,
we would be dealing not only with trees,
but also with groves (disjoint unions of trees).

The example above demonstrates how to pass
from a $k$-chain~$`g$ of partitions of~$[n]$
to a $k$-leveled tree $(T_`g,`y_`g)$ with $n$~marked leaves;
this is actually a bijection
when restricted to such chains~$`g$ that $\top \not\in `g$.
There is a unique (shortest) path from the root of $(T,`y)$ to each leaf,
and each pair of such separate at a vertex on certain level~$j$.
The labels of the two leaves
will be in the same block in the partitions~$`p_i$ for~$i\<j$,
and they will be in different blocks in~$`p_i$ for~$i>j$.
This defines the $k$-chain $`g(T,`y)$.

It will also be useful to associate with a $k$-leveled tree $(T,`y)$
a sequence $\{`l_i(T,`y)\} = \{`l_i(`g)\}$
of integer partitions as follows.
While $`l_0$ has just one part,
equal to the valency of the root of~$T$,
the partition~$`l_i$, $1 \< i \< k$,
is to have as many parts as there are vertices of $(T,`y)$ on level~$i$,
and each part is to be one less than
the number of direct descendants of the corresponding vertex.
With that,
$`r(`p_1) = `l_0(`g)$ and
$[`p_i,`p_{i+1}] \simeq L_{`l_i(`g)}$ for $1 \< i \< k$,
where $`g = `g(T,`y) = \{`p_1,\dots,`p_k\}$ and $`p_{k+1} = \top$.
For the example above,
$`l_0 = (3)$,
$`l_1 = (2)$,
$`l_2 = (1,1)$
and $`l_3 = (1,1)$.

\section{Polyscreens and colored screens}
\label{sec-screens}

\def\screen{\mathbf S}

Partition chains and leveled trees of the previous section
play in the polydiagonal compactification $\Xn$ the same role as
nests of subsets of~$[n]$ and usual trees (groves)
do in the Fulton--MacPherson compactification $X[n]$.
They index the strata and are an integral part
of the geometric description of points in $\Xn$,
explained in this section without any proofs.
It is implied by the technical work of Section~\ref{sec-strata}.

For a chain $`g = \{`p_1,\dots,`p_k\}$,
each point~$\x$ in the stratum~$S_`g$ of~$\Xn$ is represented by
a configuration~$\x'$ of distinct points in~$X$
labeled by the blocks of~$\1$
and a coherent sequence of polyscreens $\PS^{`p_1},\dots,\PS^{`p_k}$ at~$\x'$.
Let $p(\x',`b)$ be the point in the configuration~$\x'$
labeled by the block of~$\1$ that contains $`b \subseteq [n]$.
This makes sense for every block~$`b$ of every $`p \> `p_1$.

\begin{defn}
A \newterm{polyscreen} $\PS^{`p}$ at~$\x'$ is given by:
for each block $`b_i$ of~$`p$,
a configuration~$\screen_i$
of $\mbox{card}(`b_i)$ points
in the tangent space to~$X$ at $p(\x',`b_i)$,
labeled by the elements of~$`b_i$,
and
a nonzero scalar~$`a_i$,
called the \newterm{scale factor} of~$\screen_i$.
The data is considered modulo the following relations:
\begin{enumerate}
\item
translation of any screen~$\screen_i$;
\item
dilation of any screen~$\screen_i$
with compensating change of its scale factor:
$(\screen_i,`a_i)\sim(`f\screen_i,`f^{-1}`a_i)$, $`f\in\kt$;
\item
simultaneous multiplication of all~$`a_i$ by an element of~$\kt$
(rescaling).
\end{enumerate}
A sequence of polyscreens $\PS^{`p_1},\dots,\PS^{`p_k}$
is \newterm{coherent} if,
for all $j = 1,\dots,k-1$,
two labeled points in $\PS^{`p_j}$ coincide
if and only if
their labels belong to the same block of $`p_{j+1}$,
and all labeled points in $\PS^{`p_k}$ are distinct.
\end{defn}

Coherence makes a sequence $\PS^{`g}$
conform to the leveled tree $(T_`g,`y_`g)$,
as the example in Figure~8 of a point in $X\pd<9>$
does to the (right) tree in Figure~7.
This means that the root of the tree corresponds to~$X$,
each internal vertex has a screen attached to it
and the direct descendants of each vertex form a configuration
of distinct points in~$X$ or in the respective screen.
The screens in $\PS^{`g}$
attached to the phantom vertices of $(T_`g,`y_`g)$
contain just one distinct labeled point and are called \newterm{trivial};
they carry no information and are left out of the pictures.

\begin{figure}
\label{fig-polyscreens}
\xy
	(-40,0)*{\ };
	(-20,0)*{\mycaptionempty};
	(-5,29); (48,29)**\crv{(-15,18)&(4,20)&(10,26)&(55,18)&(46,26)};
	(-2,25)*{X};
	(12,26)="12357" *{\mydot};
	(23,28)		*{\mydot} +(1.3,0)*{\scriptscriptstyle 9};
	(40,24)="468"	*{\mydot};
	(10,0)="screenwd",	(0,10)="screenht",	(5,5)="gocenter",
	(7.4,-1.6)="golabel",
	(8,12)="ld",
	"ld"+"screenht"="lu",
	"ld"+"screenwd"="rd",
	"ld"+"screenht"+"screenwd"="ru",
	"ld"+"gocenter"="center",
    @i;	@={"ld","lu","ru","rd"}; s0="prev" @@{;"prev";**@{-}="prev"};
	"lu";	"12357"+(-0.5,0.1)**\dir{.};
	"ru"; 	"12357"+(0.3,0.2)**\dir{.};
	"center"+(2,3.5)*{\mydot} +(1.3,0)*{\scriptscriptstyle 7};
	"center"-(2.5,1)="15"*{\mydot};
	"center"-(0,3)="23"*{\mydot};
	(0,0)="ld",
	"ld"+"screenht"="lu",
	"ld"+"screenwd"="rd",
	"ld"+"screenht"+"screenwd"="ru",
	"ld"+"gocenter"="center",
    @i;	@={"ld","lu","ru","rd"}; s0="prev" @@{;"prev";**@{-}="prev"};
	"lu";	"15"+(-0.3,0.4)**\dir{.};
	"rd"; 	"15"+(0.4,0)**\dir{.};
	"center"-(3,-1)*{\mydot} -(0,1.3)*{\scriptscriptstyle 1};
	"center"+(3,-1)*{\mydot} -(0,1.3)*{\scriptscriptstyle 5};
	"ld"+"golabel"*{\scriptstyle\alpha=1};
	(14,-12)="ld",
	"ld"+"screenht"="lu",
	"ld"+"screenwd"="rd",
	"ld"+"screenht"+"screenwd"="ru",
	"ld"+"gocenter"="center",
    @i;	@={"ld","lu","ru","rd"}; s0="prev" @@{;"prev";**@{-}="prev"};
	"ld";	"23"+(-0.4,0)**\dir{.};
	"ru"; 	"23"+(0.2,0.4)**\dir{.};
	"center"-(3,2)*{\mydot} -(0,1.3)*{\scriptscriptstyle 2};
	"center"+(3,2)*{\mydot} -(0,1.3)*{\scriptscriptstyle 3};
	"ld"+"golabel"*{\scriptstyle\alpha=4};
	(36,0)="ld",
	"ld"+"screenht"="lu",
	"ld"+"screenwd"="rd",
	"ld"+"screenht"+"screenwd"="ru",
	"ld"+"gocenter"="center",
    @i;	@={"ld","lu","ru","rd"}; s0="prev" @@{;"prev";**@{-}="prev"};
	"lu";	"468"+(-0.5,0.1)**\dir{.};
	"ru"; 	"468"+(0.3,0.2)**\dir{.};
	"center"+(3,3)*{\mydot} +(1,-0.6)*{\scriptscriptstyle 8};
	"center"-(2,2)="46"*{\mydot};
	"ld"+"golabel"*{\scriptstyle\alpha=2};
	(28,-12)="ld",
	"ld"+"screenht"="lu",
	"ld"+"screenwd"="rd",
	"ld"+"screenht"+"screenwd"="ru",
	"ld"+"gocenter"="center",
    @i;	@={"ld","lu","ru","rd"}; s0="prev" @@{;"prev";**@{-}="prev"};
	"lu";	"46"+(-0.5,0.1)**\dir{.};
	"rd"; 	"46"+(0.3,0.2)**\dir{.};
	"center"-(3,1)*{\mydot} +(0.9,-1)*{\scriptscriptstyle 6};
	"center"+(3,1)*{\mydot} +(0.8,-0.8)*{\scriptscriptstyle 4};
	"ld"+"golabel"*{\scriptstyle\alpha=1};
	(65,27)*{\txt{Level 0}};
	(65,16)*{\txt{Level 1}};
	(65, 5)*{\txt{Level 2}};
	(65,-6)*{\txt{Level 3}};
\endxy
\end{figure}

Nontrivial screens in a polyscreen $\PS^{`p_j}$
are exactly Fulton--MacPherson screens
for those blocks of~$`p_j$ that are subdivided in $`p_{j+1}$
(all essential blocks of~$`p_j$ if $j=k$).
The data of $\PS^{`p_j}$ is equivalent to this collection of screens
together with the point in the projective space $\P^{r_j-1}$
given by the $r_j$-tuple of scale factors,
where $r_j$ is the number of nontrivial screens in~$\PS^{`p_j}$.


If $`g = \{`p_1,\dots,`p_k\}$ starts with the bottom partition $\bot$,
then for all points~$\x$ in~$S_`g$
the configuration~$\x'$ is a single point~$p$ in~$X$,
and all screens in $\PS^`g(\x)$ are based on
the {\em same} tangent space $T_p X$.
Under the additional assumption that $\operatorname{char}\Bbbk = 0$,
now made for the rest of this section,
the data of each polyscreen  $\PS^{`p_j}(\x)$ then
fits into a single {\em colored screen} $\CS^{`p_j}(\x)$.

\begin{defn}
Let a \newterm{color} be any nonempty subset of~$[n]$.
A \newterm{colored screen} $\CS^`p$ at~$p$
is a configuration of $n$~colored points $x_1,\dots,x_n$ in~$T_p X$,
considered modulo dilations of~$T_p X$,
where the color of~$x_i$ is the block of~$`p$ that contains~$i$,
such that the points of each color are centered around the origin
(their vector sum is~$0$).

A sequence $\CS^{`p_1},\dots,\CS^{`p_k}$
is \newterm{coherent} if,
for all $j = 1,\dots,k-1$,
two points of the same color coincide in $\CS^{`p_j}$ if and only if
they have the same color in $\CS^{`p_{j+1}}$,
and in  $\CS^{`p_k}$ no points of the same color coincide.
\end{defn}

To convert a polyscreen $\PS^`p$ into a colored screen $\CS^`p$,
first translate the representative screens of $\PS^`p$
to center the points around the origin,
then dilate them to make all scale factors equal.
Identifying now the underlying spaces of the screens,
place several configurations in the same~$T_p X$.
To tell them apart,
colors of points are added as a way of recording
which one of the screens each point comes from.
Figure~9 shows the simplest nontrivial example.

\begin{figure}[t]
\label{fig-conversion}
\xy
	(-13,0)*{\ };
	(75,-8)*{\mycaption{4in}{Conversion to color}};
	(14,0)="screenwd",	(0,14)="screenht",	(7,7)="gocenter",
	(10,-1.8)="golabel",
	(24,0)="stepright",
	(0,0)="ld",
	"ld"+"screenht"="lu",
	"ld"+"screenwd"="rd",
	"ld"+"screenht"+"screenwd"="ru",
	"ld"+"gocenter"="center",
    @i;	@={"ld","lu","ru","rd"}; s0="prev" @@{;"prev";**@{-}="prev"};
	"center"-(2,-2)*{\mydot} -(0,1.3)*{\scriptscriptstyle 1};
	"center"+(2,-2)*{\mydot} -(0,1.3)*{\scriptscriptstyle 2};
	"ld"+"golabel"*{\scriptstyle\alpha_{12}=2};
	(3,0)+"rd"="ld",
	"ld"+"screenht"="lu",
	"ld"+"screenwd"="rd",
	"ld"+"screenht"+"screenwd"="ru",
	"ld"+"gocenter"="center",
    @i;	@={"ld","lu","ru","rd"}; s0="prev" @@{;"prev";**@{-}="prev"};
	"center"-(3,1.5)*{\mydot} -(0,1.3)*{\scriptscriptstyle 3};
	"center"+(3,1.5)*{\mydot} -(0,1.3)*{\scriptscriptstyle 4};
	"ld"+"golabel"-(1,0)*{\scriptstyle\alpha_{34}=-1};
	"rd"+(5,6)*{\mapsto};
	"stepright"+"ld"="ld",
	"ld"+"screenht"="lu",
	"ld"+"screenwd"="rd",
	"ld"+"screenht"+"screenwd"="ru",
	"ld"+"gocenter"="center",
    @i;	@={"ld","lu","ru","rd"}; s0="prev" @@{;"prev";**@{-}="prev"};
	"center"-(4,-4)*{\mydot} -(0,1.3)*{\scriptscriptstyle 1};
	"center"+(4,-4)*{\mydot} -(0,1.3)*{\scriptscriptstyle 2};
	"ld"+"golabel"*{\scriptstyle\alpha_{12}=1};
	(3,0)+"rd"="ld",
	"ld"+"screenht"="lu",
	"ld"+"screenwd"="rd",
	"ld"+"screenht"+"screenwd"="ru",
	"ld"+"gocenter"="center",
    @i;	@={"ld","lu","ru","rd"}; s0="prev" @@{;"prev";**@{-}="prev"};
	"center"+(3,1.5)*{\mydot} -(0,1.3)*{\scriptscriptstyle 3};
	"center"-(3,1.5)*{\mydot} -(0,1.3)*{\scriptscriptstyle 4};
	"ld"+"golabel"*{\scriptstyle\alpha_{34}=1};
	"rd"+(5,6)*{\mapsto};
	"stepright"+"ld"="ld",
	"ld"+"screenht"="lu",
	"ld"+"screenwd"="rd",
	"ld"+"screenht"+"screenwd"="ru",
	"ld"+"gocenter"="center",
    @i;	@={"ld","lu","ru","rd"}; s0="prev" @@{;"prev";**@{-}="prev"};
	"center"-(4,-4)*{\mydot} -(0,1.3)*{\scriptscriptstyle 1};
	"center"+(4,-4)*{\mydot} -(0,1.3)*{\scriptscriptstyle 2};
	"center"+(3,1.5)*{\scriptstyle\circ} -(0,1.3)*{\scriptscriptstyle 3};
	"center"-(3,1.5)*{\scriptstyle\circ} -(0,1.3)*{\scriptscriptstyle 4};
\endxy
\end{figure}

Since this conversion of polyscreens into colored screens respects coherence,
points in~$\Xn$ corresponding to collisions at a single point in~$X$
can be viewed in terms of coherent sequences of colored screens.
This interpretation is useful in Section~\ref{sec-isotropy}
for studying the natural action of~$\Sym_n$ on~$\Xn$.

\section{Construction of the compactification}
\label{sec-construct}

For a partition~$`p$ of~$[n]$,
denote by $`D^`p\subseteq X^n$ the subset of all points $(x_1,\dots,x_n)$
with $x_i=x_j$ whenever $i$ and~$j$ are in the same block of~$`p$,
and call $`D^`p$ a \newterm{ polydiagonal}.
The diagonals of~$X^n$ correspond to
partitions with only one essential block.
The set of all polydiagonals in~$X^n$ is naturally
a lattice isomorphic to $\Ln$,
with its top element~$X^n$ itself.

\begin{thm}
\label{construction}
The following\/ $(n-1)$-stage sequence of blowups
results in a smooth compactification $\Xn$ of
the configuration space of\/~$n$ distinct labeled points
in a~smooth algebraic variety~$X$:
\begin{itemize}
\item
the first stage is the blowup of\/~$`D$, the small diagonal of\/~$X^n$;
\item
the $k$-th stage, $1<k<n$,
is the blowup of the disjoint union of
the previous stage proper transforms $Y^`p_{k-1}$ of\/ $`D^`p$,
for all partitions~$`p$ of the set\/ $[n]=\{1,\dots,n\}$
into exactly\/~$k$ blocks.
\end{itemize}
\end{thm}

\begin{rem}
In the language of De~Concini and Procesi~\cite{DeCP},
the building set for this iterated blowup construction consists of
{\em all} possible intersections of the diagonals of~$X^n$,
and therefore it is maximal.
The building set of the Fulton--MacPherson compactification includes only
those intersections that fail to be normal,
so $X[n]$ is the minimal compactification of $\Conf(X,n)$
with the property that the complement to the configuration space
is a divisor with normal crossings.
\end{rem}

Two smooth subvarieties $U$~and~$V$ of a smooth algebraic variety~$W$
are said to \newterm{intersect cleanly}
if $U\not\subset V\not\subset U$,
their scheme-theoretic intersection is smooth
and the tangent bundles satisfy $T(U\cap V)=TU\cap TV$.
Two polydiagonals $`D^\1$ and $`D^\2$ in~$X^n$ intersect cleanly
unless one of them contains the other;
the noncontainment condition is that
the partitions $\1$~and~$\2$ are incomparable in~$\Ln$.

Recall two standard results about
the behaviour of clean intersections under blowups:

\begin{lem}
\label{clean}
Let\/~$W$ be a smooth algebraic variety and
let\/ $U$,~$V$ be smooth subvarieties of\/~$W$ intersecting cleanly.
Then
\begin{enumerate}
\item
\label{be-disjoint}
the proper transforms of\/ $U$ and\/~$V$ in $\Bl_{U\cap V}W$ are disjoint;
\item
\label{ZinUcapV}
if\/~$Z$ is a smooth subvariety of~$U\cap V$,
then the proper transforms of\/ $U$ and\/~$V$ in $\Bl_Z W$
intersect cleanly.\qed
\end{enumerate}
\end{lem}

\begin{proof}[Proof of Theorem 1.]
Denote the space obtained at stage~$k$ by~$Y_k$
and organize the projections of the fiber squares of all stages as
\begin{equation}
\label{sequence}
\xymatrix @-3pt
{
\Xn=
Y_{n-1}\ar[r] &
Y_{n-2}\ar[r] &
{}\dots\ar[r] &
Y_1    \ar[r] &
Y_0=X^n.
}
\end{equation}
Then $Y^`p_0=`D^`p$ and $Y^`p_k$ is
the proper transform of~$Y^`p_{k-1}$ in~$Y_k^{\vphantom{`p}}$
if $`r(`p)\ne k$,
while $Y^`p_{`r(`p)}$ is
the component of the exceptional divisor over $Y^`p_{`r(`p)-1}$.

The statement will follow once it has been shown that
the stated sequence of blowups can indeed be performed.
For this,
it suffices to check that the centers of
those simultaneous blowups will have indeed become disjoint
after the previous stages of the construction.
The proof will be done by induction on~$k$,
for all $\Xn$ at the same time;
after stage~$k$, the induction will stop for $X\pd<k+1>$,
and it will continue on for $\Xn$ with $n>k+1$.

For any pair of distinct partitions $\1$ and~$\2$ of~$[n]$ into two blocks,
their meet $\1\wedge\2$ is the \inquotes{nonpartition}\negthickspace,
so $`D^\1\cap`D^\2=`D^{\1\wedge\2}=`D$,
the small diagonal of~$X^n$.
By Lemma~\ref{clean}\ref{be-disjoint},
the transforms $Y_1^\1$ and $Y_1^\2$ will be disjoint,
making the second stage possible.

Assume that stage $k-1$ has been performed;
this means that the varieties $\Xn$ have been constructed
for $1\leqslant n\leqslant k$,
and only those for $n>k$ are still being built.
Also assume that the proper transforms
$Y^`p_{k-1}$ for~$`p$ with $`r(`p)=k$ are disjoint.

For each partition $`p\in \Ln$ with $`r(`p)=k$,
the projection $X\pd<k>\to X^k$ pulls back
the obvious isomorphism $X^k \simeq `D^`p \subset X^n$
to an isomorphism $X\pd<k> \simeq Y_{k-1}^`p \subset Y_{k-1}^{\vphantom{`p}}$.
All these subvarieties are disjoint by the inductive assumption,
and can all be blown up at the same time.
This defines the variety $X\pd<k+1>$.

To provide the inductive step necessary to continue
the construction of $\Xn$ for $n>k+1$,
the intersection $Y_k^\1\cap Y_k^\2$ must be empty
for all pairs of distinct $\1$, $\2$ in $\Ln$
with $`r(\1)=`r(\2)=k+1$.
Such a pair automatically satisfies the noncontainment condition;
so $`D^\1\cap`D^\2=`D^{\1\wedge\2}$ is a clean intersection.
Since $`r=`r(\1\wedge\2)<k+1$,
a repeated use of Lemma~\ref{clean}\ref{ZinUcapV}
shows that $Y^\1_{`r-1}\cap Y^\2_{`r-1}=Y^{\1\wedge\2}_{`r-1}$ is
a clean intersection,
and then Lemma~\ref{clean}\ref{be-disjoint}
implies that $Y^\1_`r\cap Y^\2_`r$ is empty.
The proper transforms of $`D^\1$ and $`D^\2$
become disjoint after stage $`r\<k$,
and the proof is complete.
\end{proof}

\begin{corol}
\label{middleX}
For each $`p\in\Ln$ we have $Y_{`r(`p)-1}^`p\simeq X\pd<`r(`p)>$.
\end{corol}

\begin{proof}
This has been obtained while proving the theorem,
and is formulated separately only for the ease of future reference.
\end{proof}

\begin{FBL}
Let\/ $V^1_0\subset V^2_0\subset\dots\subset V^s_0\subset W_0^{\vphantom{s}}$
be a flag of smooth subvarieties in a smooth algebraic variety~$W_0$.
For $k=1,\dots,s$, define inductively:
$W_k$ as the blowup of\/ $W_{k-1}$ along $V^k_{k-1}$;
$V^k_k$ as the exceptional divisor in~$W_k$;
and\/ $V^i_k$, for $i\ne k$,
as the proper transform of\/~$V^i_{k-1}$ in~$W_k$.
Then the preimage of\/~$V^s_0$ in the resulting variety~$W_s$
is a normal crossing divisor $V_s^1\cup\dots\cup V_s^s$.
\end{FBL}

\begin{rem}
This auxiliary result is implicit in earlier works
\cite{FM,Kapranov:Chow}.
\end{rem}

\begin{proof}
In a blowup $p \colon \Bl_Z W \to W$
of a smooth algebraic variety~$W$ along a smooth center~$Z$,
if~$\tilde{V}$ is the proper transform of a smooth variety~$V \supset Z$,
then in terms of ideal sheaves $\I(p^{-1}(V))=\I(\tilde{V})\cdot\I(E)$.
Applied at each step,
this equality yields
$\I(p_s^{-1}(V^s_0)) = \I(V^1_s) \times\dots\times \I(V^s_s)$,
where $p_s\colon W_s\to W_0$ denotes the composition of the stated blowups.
\end{proof}

\begin{prop}
\label{divisor}
For each partition~$`p$ of\/~$[n]$ with at least one essential block,
there is a smooth divisor $D^`p\subset \Xn$ such that:
\begin{enumerate}
\item
The union of these divisors is $D=\Xn\smallsetminus \Conf(X,n)$.
\item
Any set of these divisors meets transversally.
\item
An intersection $D^{`p_1}\cap\dots\cap D^{`p_k}$ of divisors
is nonempty exactly when the partitions form a chain.
In other words,
the incidence graph of\/~$D$
coincides with the comparability graph of
the lattice~$\Ln$ with the top partition removed.
\end{enumerate}
\end{prop}

\begin{corol}
\label{strata}
\begin{enumerate}
\item
$\Xn$ is stratified by strata $\,S_`g=\bigcap_{`p\in`g}D^`p$
parametrized by all chains~$`g$ in~$\Ln$.
\item
The codimension of\/~$S_`g$ in~$\Xn$
is equal to the length of\/~$`g$.
\item
The intersection of two strata $S_`g$~and\/~$S_{`g'}$
is nonempty exactly when $`g\cup`g'$ is a chain,
in which case $S_`g\cap S_{`g'}=S_{`g\cup`g'}$.
In particular,
$S_`g\supset S_{`g'}$ if and only if\/~$`g\subset`g'$.
\end{enumerate}
\end{corol}

\begin{proof}
We concentrate on the normal crossing property,
which implies the other claims.

By construction,
the proper transform of every polydiagonal $`D^`p\subset X^n$
under $\Xn\to X^n$ is a smooth divisor;
it will be denoted by~$D^`p$.
The proper transforms of $`D^\1$ and $`D^\2$ become disjoint
when that of their intersection $`D^{\1\wedge\2}$ is blown up,
unless one of $`D^\1$ and $`D^\2$ contains the other,
that is,
unless $\{\1,\2\}$ is a chain.

In order to show that for any saturated (maximal length) chain $`g=\{`p_i\}$,
the union $D^`g=D^{`p_1}\cup\dots\cup D^{`p_{n-1}}$
is a normal crossing divisor in~$\Xn$,
consider the flag of polydiagonals
$`D^{\1}\subset\dots\subset`D^{`p_{n-1}}\subset X^n$.
The blowups of $Y^`p_{`r(`p)-1}$ for $`p\not\in`g$
are irrelevant for the intersection of the components of~$D^`g$
because their centers are disjoint from
$$
\bigcap_{i=1}^{n-1} Y^{`p_i}_{`r(`p)-1};
$$
hence,
the Flag Blowup Lemma can be applied.
The normal crossing property of $D^`g$ follows by the lemma,
and so does the proposition:
since any chain~$`g'$ is refined by a saturated chain~$`g$,
the components of $\bigcup_{`p\in`g'} D^`p$
form a subset of components of $\bigcup_{`p\in`g} D^`p$.
\end{proof}

\subsection*{Enumeration of the strata.}
The number of strata in~$\Xn$, $n>1$,
is equal to the number $2Z(n)$ of chains in~$\Ln$.
There is a factor of $2$ here
because half of the chains contain $\bot$ and half do not
(the top~$\top$ is always excluded).
Sloane and Plouffe~\cite{EIS} catalogued
the sequence $\{Z(n)\}$ of integers as M3649.
Since the following recurrence relation is immediate:
$$
\label{recurrenceZn}
Z(n)=\sum_{k=1}^{n-1}S(n,k)Z(k),
$$
the first few values of $Z(n)$ are easy to compute.
No closed general formula is known,
although Babai and Lengyel described the asymptotics of $Z(n)$,
up to yet undetermined constant \cite{BabaiLengyel,Lengyel}.

Here is a small table of the numbers of strata in $X[n]$ and $\Xn$:
$$
\begin{tabular}{|c|c|c|c|c|c|c|c|c|}
\hline
$n$ &2 &3 &4 &5 &6 &7 &8 &9 \\
\hline
$X[n]$ &2 &8 &52 &472 &5504 &78416 &1320064 &25637824 \\
\hline
$\Xn$ &2 &8 &64 &872 &18024 &525520 &20541392 &1036555120 \\
\hline
\end{tabular}
$$

As codimension-$1$ strata
are the components~$D^`p$ of the divisor at infinity,
there are $B(n)-1$ of them,
where $B(n)$ is the Bell number,
equal to the number of partitions of~$[n]$.
The minimal strata have codimension $n-1$ and
correspond to saturated chains in~$\Ln$,
whose number is $2^{1-n}n!(n-1)!$.

\section{The \Hopo\ of $\Xn$}
\label{sec-hodge}

If~$X$ is a smooth complex algebraic variety,
the construction of~$\Xn$ allows an easy derivation
of a formula for the \Hopo,
hence, for the \Popo\ of~$\Xn$ in terms of those of~$X$.

The notion of a {\em virtual \Popo\/} extends
the usual one to all complex algebraic varieties and
provides a good tool for computing the \Popo s of blowup constructions.

\begin{lem}
\label{popo}
\begin{enumerate}
\item
If\/~$Y$ is smooth and compact,
the virtual \Popo\ $P(Y)$ coincides with the usual \Popo\ of\/~$Y$.
\item
If\/ $Z$ is a closed subvariety of\/ $Y$,
then $P(Y)=P(Z)+P(Y\smallsetminus Z)$.
\item
If\/ $Y'\to Y$ is a bundle with fiber $F$
which is locally trivial in the Zariski topology,
then $P(Y')=P(Y)P(F)$.
\qed
\end{enumerate}
\end{lem}

Using Deligne's mixed Hodge theory \cite{Deligne:Hodge,Deligne:poids},
Danilov and Khovanskii defined a refinement of~$P(X)$,
the {\em virtual \Hopo\/} $e(X)$,
also called the {\em Serre polynomial},
and proved~\cite{DanilovKhovanskii}
that it has the properties listed in Lemma~\ref{popo},
also independently found by Durfee~\cite{Durfee}.
Cheah, Getzler and Manin computed
the \Hopo s of the Fulton--MacPherson compactifications
via generating functions \cite{Cheah,Getzler,Manin},
while the original paper dealt with
summation over trees (groves).

\begin{prop}
For any two positive integers $m$~and\/~$n$,
there is a polynomial\/ $U^m_n(t,x)$ such that
for any smooth $m$-dimensional complex algebraic variety~$X$
the \Hopo\ of $\Xn$ is $e(\Xn;z,\bar{z})=U^m_n(z\bar{z},e(X;z,\bar{z}))$,
and in particular,
$P(\Xn;t)=U^m_n(t,P(X;t))$.
The polynomials $U^m_n(t,x)$ satisfy the recurrence relation
$$
U^m_n(t,x) = x^n + \sum_{k=1}^{n-1} S(n,k) h_{(n-k)m}(t) U^m_k(t,x),
$$
where $h_d(t)=P(\C\P^{d-1})-1=t^{2d-2}+\dots+t^4+t^2$.
\end{prop}

\begin{proof}
Straightforwardly from the construction of~$\Xn$ and Lemma~\ref{popo},
$$
e(Y_k) =
e(Y_{k-1}) + \sum_{`r(`p)=k} \big(e(\P(N^`p))-1\big) e(Y^`p_{k-1}),
$$
where~$N^`p$ is the fiber of
the normal bundle to $Y^`p_{k-1}$ in~$Y^{\phantom{`p}}_{k-1}$,
which is $\C^{(n-k)m}$ by an easy dimension count.
Corollary~\ref{middleX} converts this formula into
$$
e(Y_k) = e(Y_{k-1}) + S(n,k) h_{(n-k)m}(z\bar{z}) e(X\pd<k>).
$$
Since $Y_0=X^n$ and $Y_{n-1}=\Xn$,
there results a recurrence relation
$$
e(\Xn)=e(X)^n + \sum_{k=1}^{n-1} S(n,k) h_{(n-k)m}(z\bar{z}) e(X\pd<k>),
$$
and both claims immediately follow.
\end{proof}

A nonrecursive expression for $U^m_n(t,x)$ can be found
by expanding in the right-hand side of the recurrence
the terms with the highest~$k$ present in a loop down to $k=2$ terms:
$$
U^m_n(t,x)= x^n+
\sum_{s=1}^{n-1}
    \Bigg(\!
	x^s
	\sum_{r=1}^{n-s}
	    \sum_{\mathrm J_{s,n}^{r}}
		\prod_{i=1}^{r}
		    S(j_i,j_{i-1}) h_{(j_{i}-j_{i-1})m}(t)
    \!\Bigg),
$$
where $\mathrm J_{s,n}^{r}=\{(j_0,\dots,j_r)\in`Z^{r+1}
\suchthat s=j_0<\dots< j_r=n\}$.

Similar computations of the \Hopo s of the strata in
the stratification of~$\Xn$ from Corollary~\ref{strata}
can be carried out using the description of their structure
given in Section~\ref{sec-strata}.

\section{$\Xn$ as a closure and a surjection $\Xn\to X[n]$}
\label{sec-closure}

In this section I present~$\Xn$ as the closure
of the configuration space embedded in a product of blowups,
exhibit a surjection $\Xn\to X[n]$,
and write it as an iterated blowup.

First, the results of Section~\ref{sec-construct}
about the structure of~$\Xn$ at infinity
should be rephrased in terms of ideal sheaves.
Let $\I(`D^`p)$ be the ideal sheaf of
$`D^`p$ in~$\O_{X^n}$.
For any~$k$, $1\leqslant k\leqslant n-1$,
let $`t_k\colon Y_k\to X^n$ be the appropriate composition of projections
from Eq.~(\ref{sequence}),
let $\I_k(`p)$ be the ideal sheaf in $\O_{Y_k}$
generated by $`t^*_k(\I(`D^`p))$,
and also let $\I(Y^`p_k)$ be the ideal sheaf of $Y^`p_k$ in~$\O_{Y_k}$.
This notation,
although similar to \FM's,
is not quite the same.
The assertions of Proposition~\ref{divisor} can be restated as
$$
\I_{n-1}(`p)=\prod_{`p'\leqslant`p}\I(D^{`p'}),
$$
while at the intermediate stages
$$
\I_{k}(`p)=\prod_{`p'\leqslant`p\text{ with }`r(`p')\leqslant k}
\I(Y^{`p'}_k).
$$
Since $Y^{`p'}_k\subset Y^{\vphantom{`p'}}_k$ is
a divisor if $`r(`p')<k$,
it follows that $\I_{k}(`p)=\I(Y^{`p}_k)\cdot\mathcal J$,
where~$\mathcal J$ is an invertible ideal sheaf.

\begin{prop}
\label{closure}
The variety~$\Xn$ constructed by blowing up is 
the closure of the configuration space $\Conf(X,n)$ in\/
$$
\prod_{`p\in\Ln}\Bl_{`D^`p}X^n.
$$
\end{prop}

\begin{rem}
The top partition contributes the factor $X^n$.
\end{rem}

\begin{proof}
By induction on~$k$,
each $Y_k$ is the closure of~$\Conf(X,n)$ in
$$
X^n\times \prod_{`r(`p)\<k}\Bl_{`D^`p}X^n.
$$
The basis is clear: $Y_0=X^n$.
Then,~$Y_k$ is the blowup of  $Y_{k-1}$ along
$$
\coprod_{`r(`p)=k} Y^`p_{k-1},
$$
or in other terms,
along
$$
\I\Bigg(\coprod_{`r(`p)=k}Y^`p_{k-1}\Bigg)=\prod_{`r(`p)=k}\I(Y^`p_{k-1}).
$$
This ideal sheaf becomes
$$
\I_{k-1}=\prod_{`r(`p)=k}\I_{k-1}(`p)
$$
upon multiplying by an invertible ideal sheaf,
and blowing up $\I_{k-1}$ is equivalent to
taking the closure of the graph of the rational map from $Y_{k-1}$ to
$$
\prod_{`r(`p)=k}\Bl_{`D^`p}X^n.
$$
This provides the inductive step,
and eventually $Y_{n-1}=\Xn$.
\end{proof}

Both the statement and its proof parallel those by \FM\ \cite[Prop.~4.1]{FM},
who use pullbacks by $X[n]\to X^n\to X^S$,
for $S\subset[n]$, $\#S>1$,
and also by $f_S\colon Y_k\to X^n\to X^S$ at the intermediate stages,
while here $`t_k\colon Y_k\to X^n$.
A slight reformulation of their characterization of~$X[n]$ as a closure
elucidates its similarity with $\Xn$.
For each~$S$ as before,
take the diagonal $`D^S\subset X^n$ and pull back its ideal sheaf
by the first of the two arrows whose composition is~$f_S$;
this gives the same ideal sheaf $f^*_S(\I(`D))$.

\begin{prop}
\label{reformulated}
The variety $X[n]$ is 
the closure of $\Conf(X,n)$ in
$$
X^n\times\prod_{S\subset[n],\#S>1}\Bl_{`D^S}X^n.
$$
\end{prop}

The two compactifications can now be related.

\begin{prop}
\label{ontoFM}
For each~$n \> 1$, there is a surjection $\vartheta_n\colon \Xn\to X[n]$.
\end{prop}

\begin{proof}
Start with notation for the products
from Propositions \ref{closure}~and~\ref{reformulated}:
$$
\Pi = \prod_{`e(`p)\geqslant 1} \Bl_{`D^`p}X^n,
\quad \text{ and }\quad
\Pi' = \prod_{`e(`p)=1} \Bl_{`D^`p}X^n,
$$
where $`e(`p)$ is
the number of essential blocks in a partition~$`p$.
If~$S$ is the only essential block of~$`p$,
then $`D^S=`D^`p$,
so $\Pi'$ can indeed be used for $X[n]$.
\begin{equation*}
\mbox{
\xymatrix @-4mm{
    &
    &
      {{}\Pi}
      \ar^p[dd]
    \\
      {}\Conf(X,n)
      \ar[r]
      \ar@/^/[urr]
      \ar@/_/[drr]
    &
      X^n
      \ar@{.>}_{\!\!\!\!\phi}[ur]
      \ar@{.>}^{\!\!\!\!\psi}[dr]
    &
    \\
    &
    &
      {{}\Pi'}
  }
}
\phantom{mmmmm}
\mbox{
\xymatrix @-4mm{
    &
    &
      X^n\times \Pi
      \ar^{\text{id}\times p}[dd]
    \\
      {}\Conf(X,n)
      \ar[r]
      \ar@/^.7pc/[urr]
      \ar@/_.7pc/[drr]
    &
      X^n
      \ar@{.>}_{\!\!\!G(\phi)}[ur]
      \ar@{.>}^{\!\!\!G(\psi)}[dr]
    &
    \\
    &
    &
      X^n\times \Pi'
  }
}
\end{equation*}

Now take the left of these two diagrams,
where $`f$ and~$`q$ are rational maps defined on $\Conf(X,n)$,
and notice that the projection $\mbox{id}\times p$ maps
the closure $\overline{G(`f)}$ of the graph of~$`f$ onto
the closure $\overline{G(`q)}$ of the graph of~$`q$.
\end{proof}

This surjection $\vartheta_n$ admits a more explicit description.
For~$n\leqslant 3$,
it is the identity map;
otherwise,
it can be written as a composition
\begin{equation}
\label{theta}
\xymatrix @+3pt
{
X\pd<n>=
W_{n-2}\ar[r]^<(0.3){\beta_{n-2}} &
W_{n-3}\ar[r]^<(0.35){\beta_{n-3}} &
{}\dots\ar[r]^{\beta_{3}} &
W_2    \ar[r]^<(0.25){\beta_{2}} &
W_1=X[n],
}
\end{equation}
where $W_k \xrightarrow{`b_k} W_{k-1}$ is
the blowup in~$W_{k-1}$ of (the disjoint union of
the proper transforms under $`b_{k-1}\circ\dots\circ `b_2$ of)
some strata $X(\S)$ of~$X[n]$;
their encoding nests~$\S$ are characterized below.
Favoring imprecision over repetitiveness,
I will neglect to reiterate the ritual phrase
that in the previous sentence appears in parentheses.

Let $U\subset X[n]$ be the union of all strata $X(\S)$
such that the nest $\S$ contains two disjoint subsets of~$[n]$.
The irreducible components of this codimension~2 reduced subscheme are
$X(\S)$ for all nests $\S=\{S_1,S_2\}$ with $S_1\cap S_2=\varnothing$,
which intersect transversally \cite[Theorem~3]{FM}.
The map~$\vartheta_n$ is an iterated blowup of $X[n]$ along~$U$,
but not all the strata contained in~$U$ are centers of a blowup $`b_k$.
The components of the center of $`b_k$ are the strata $X(\S)$
such that $\S$ is the set of all essential blocks of
a partition $`p\in\Ln$ with $`r(`p)=k$ and $`e(`p)>1$,
which is always a nest.
The transversality of the strata guarantees that,
whenever the sequence~$`b_2,\dots,`b_{n-2}$ calls for
two intersecting strata to be in the center of the same~$`b_k$,
the previous stages will have made them disjoint.
The sequence itself implies that,
whenever $X(\S)\subset X(\S')$ are both to become centers,
the smaller stratum is blown up before the larger one.

Alternatively,
the variety $W_k$ can be defined as the closure of $\Conf(X,n)$ in
$$
X^n\times \prod_{`r(`p)\<k\text{ or }`e(`p)=1} \Bl_{`D^`p}X^n,
$$
and an argument similar to Proposition~\ref{closure}
shows that this is equivalent to the blowup description.

\begin{examples}
Here $X(S_1,\dots,S_k)=D(S_1)\cap\dots\cap D(S_k)$
refers to strata of~$X[n]$.

The map $\vartheta_4$ is the blowup of
3 disjoint codimension-2 strata $X(12,34)$ and alike,
for the nests obtained from the 3 partitions of shape $(2,2)$.
The divisor $D^{12|34}\subset X\pd<4>$ is
a $\P^1$-bundle over $X(12,34)$.

For $n=5$,
there are two maps in Eq.~(\ref{theta}).
The first blows up 10 disjoint codimension-2 strata,
like $X(123,45)$, corresponding to
the partitions of shape $(3,2)$.
The second blows up 15 disjoint codimension-2 strata,
like $X(12,34)$,
corresponding to $(2,2,1)$.

For $n=6$,
there are three stages according to the partitions
$$
(4,2), \ \ \ (3,3); \quad\quad\quad
(3,2,1), \ \ \ (2,2,2); \quad\quad\quad
(2,2,1,1).
$$
Here we encounter inclusions like $X(12,34,56)\subset X(12,34)$.
Interestingly,
the proper transform by $\vartheta_6$ of $X(12,34,56)$,
which is the divisor $D^{12|34|56}$,
is a bundle over $X(12,34,56)$ with fiber
$\P^2$~blown up at three points.
Proposition~\ref{small-bricks} generalizes this observation.
\end{examples}

The preimage in~$\Xn$ of a stratum of~$X[n]$ is
$$
\vartheta_n^{-1}(X(\S))=\bigcup_{(T,`y)\in`h^{-1}(T(\S))} S_{(T,`y)},
$$
the union of all strata encoded by the leveled trees $(T,`y)$
with the same base tree $T(\S)$ and
any legal assignment of levels to its interior vertices.
The map~$\vartheta_n$ is thus strata-compatible.

\begin{prop}
The fibers of $\vartheta_n \colon \Xn \to X[n]$ are independent of~$X$
and even of its dimension.
The fiber over a point in $X(\S)$ that is not in any smaller stratum
is completely determined by the nest~$\S$.
\end{prop}

\begin{proof}
The normal space~$N_\x$ at a point~$\x$ to $X(\S) \subset X[n]$
is independent of $\dim X$ (assumed positive):
its dimension is equal to the cardinality of the nest~$\S$.
The nest alone determines the iterated blowup of~$N_\x$
induced from Eq.~(\ref{theta}),
and the preimage of the origin under it
is isomorphic to $\vartheta_n^{-1}(\x)$.
\end{proof}

\section{Structure of the strata}
\label{sec-strata}

\def\1{{`p_1}}
\def\2{{`p_2}}
\def\Br{\mathcal B}

This section begins by discussing
a family of linear subspace arrangements indexed by integer partitions;
each of them leads to a projective variety that will be called a brick.
Points of a brick correspond to polyscreens,
and by presenting the strata of~$\Xn$ as bundles over~$X\pd<k>$
whose fibers are products of bricks,
the polyscreen description of~$\Xn$ is established here.

The configuration space $\Conf(\mathbb A^{\!1},n)$ is the complement to
the braid arrangement of hyperplanes in $\mathbb A^{\!n}$,
the motivating example for much of
the theory of hyperplane arrangements~\cite{OrlikTerao}.
The analogue for $\Am$,
denoted by~$\bar\Br^m_n$,
is an arrangement of codimension~$m$ linear subspaces of~$(\Am)^n$.
Its strata are various intersections of the large diagonals,
so the partitions of~$[n]$ index them,
for each $m \> 1$;
in other words,
the intersection lattice of~$\bar\Br^m_n$
is isomorphic to the partition lattice~$\Ln$.
These and all other subspace arrangements encountered in this section are
$c$-plexifications of hyperplane arrangements~\cite{Bjorner}.
This means practically that most information about~$\bar\Br^m_n$
can be extracted from
the braid arrangement~$\bar\Br^1_n$.

For any partition~$`p$ of~$[n]$,
the images in the quotient $C^m_`p = (\Am)^n / `D^`p$
of those large diagonals that contain $`D^`p$
form an induced arrangement $\Br^m_`p$.
For $`p = \bot(\Ln)$,
it is denoted by $\Br^m_{n-1}$
(actual subscripts will be integers~$`n_i$);
if~$m=1$,
this is the Coxeter arrangement of type~$A_{n-1}$.
For other partitions,
$\Br^m_`p$ is a product arrangement,
as Lemma~\ref{product-arrangement} shows below.

\def\arr{\mathcal A}
\def\8{{\vphantom{1}}}

For two subspace arrangements $\arr_i = \{K^i_1,\dots,K^i_{s_i}\}$
in $\Bbbk$-vector spaces~$V_i$,
$i=1,2$,
the product arrangement $\arr_1 \times \arr_2$ in $V_1 \oplus V_2$
is the collection of subspaces
$\{K^1_1 \oplus V^\8_2, \dots, K^1_{s_1} \oplus V^\8_2,
   K^2_1 \oplus V^\8_1, \dots, K^2_{s_2} \oplus V^\8_1\}$.
For each integer partition $`l = (`n_1,\dots,`n_r)$,
define~$\Br^m_`l$ as the product 
$\Br^m_{`n_1} \times\dots\times \Br^m_{`n_r}$.
The intersection lattice of a product is the product of those of the factors;
for~$\Br^m_`l$ this gives the lattice
$L_`l = L_{[`n_1+1]} \times\dots\times L_{[`n_r+1]}$.

As an example,
take for $`l$ the finest partition $(1,\dots,1)$ of~$r$,
often denoted by $1^r$.
Since~$\Br^1_1$ is the arrangement $\{0\}$ in~$\Bbbk$,
its $r$-th power $\Br^1_{1^r}$
is the arrangement of coordinate hyperplanes in~$\Bbbk^r$.

\begin{lem}
\label{product-arrangement}
Up to a change of coordinates
$\Br^m_`p \simeq \Br^m_`l$,
where $`l=`l(`p)$ is the essential shape of~$`p$.
\end{lem}

\begin{proof}
Look at the equations of the large diagonals containing~$`D^`p$,
that is,
$`D^{ij} = \{(x_1,\dots,x_n) \in (\Am)^n \suchthat x_i = x_j\}$
for all pairs of $i$~and~$j$ belonging to the same block of~$`p$.
Equations coming from different blocks of~$`p$ are independent of each other,
leading to the product decomposition.
\end{proof}

The polydiagonal compactification $\Am\pd<n>$
is the maximal blowup of the arrangement~$\bar\Br^m_n$,
in the sense that all strata of~$\bar\Br^m_n$
are blown up in the course of its construction.
In the same fashion,
all strata of the arrangement~$\Br^m_`l$
can be blown up in the ascending order given by their dimensions.
The first stage is always the blowup of the origin,
creating the exceptional divisor $\P(C^m_`l) \simeq \P^{m|`l|-1}$,
where $|`l|$ is the sum of all parts of~$`l$.

The main objects of interest for this section are defined as follows.

\begin{defn}
For any integer partition~$`l$,
a \newterm{brick} $M^m_`l$ is the proper transform of $\P(C^m_`l)$
in the maximal blowup of $\Br^m_`l$.
If~$`l$ has only one part,
the brick~$M^m_`l$ is \newterm{simple},
otherwise it is \newterm{compound}.
The \newterm{open} brick $\0 M^m_`l$ is
the complement in $\P(C^m_`l)$ of the projectivization of~$\Br^m_`l$.
\end{defn}

\begin{examples}
The brick $M^m_1$ is just $\P^{m-1}$ (a single point if $m=1$).

The bricks $M^m_2$ and $M^m_{1,1}$ are blowups of $\P^{2m-1}$;
their centers are,
respectively,
three and two copies of $M^m_1$.

The bricks $M^m_3$, $M^m_{2,1}$ and $M^m_{1,1,1}$ are 
2-stage blowups of $\P^{3m-1}$;
the lower intervals in $L_3$, $L_{2,1}$ and $L_{1,1,1}$
determine their centers,
respectively:
\begin{list}{}{}
\item 
7 copies of $M^m_1$,
then 6 copies of $M^m_2$;
\item 
4 copies of $M^m_1$,
then 3 copies of $M^m_{1,1}$ and 1 copy of $M^m_2$;
\item 
3 copies of $M^m_1$,
then 3 copies of $M^m_{1,1}$.
\end{list}
For $M^1_3$,
look again at Figures 2~and~3 
on page~\pageref{fig-x4}.
Similar pictures for $M^1_{2,1}$ and  $M^1_{1,1,1}$
are in Figures 10~and~11.
Comparison of these figures suggests that
refining the indexing partition corresponds to
omitting some subspaces from the arrangement.
This is proved in general in Proposition~\ref{sub-and-super}.

\begin{figure}[t]
\label{fig-four-lines} 		
\def\0{\kern-1pt}
\xy
	(16,-12)*{\mycaptionempty};
	(0,-6); (0,18)**\dir{-};
	(-14,0); (14,0)**\dir{-};
	(11.4,-5.1); (-3.6,17.4)**\dir{-};
	(-11.4,-5.1); (3.6,17.4)**\dir{-};
	(0,0)*{\mydot};
	(8,0)*{\mydot};
	(0,12)*{\mydot};
	(-8,0)*{\mydot};
	(3.9,15.1)*{\scriptstyle 1\!2};
	(9,-4)*{\scriptstyle 1\!3};
	(-12.6,1.2)*{\scriptstyle 4\05};
	(-1.4,-4.2)*{\scriptstyle 2\03};
	(24,-4)="m",
	(8,0)+"m"="bot",	
	(0,7)+"m"="12",
	(4,7)+"m"="13",
	(8,7)+"m"="23",
	(16,7)+"m"="45",
	(0,14)+"m"="123",
	(8,14)+"m"="1245",
	(12,14)+"m"="1345",
	(16,14)+"m"="2345",
	(8,21)+"m"="top",
	@={"top","123","12","bot","13","123","23","bot","45","2345","23",
		"2345","top","1245","12","1245","45","1345","13","1345","top"},
	s0="prev" @@{;"prev";**@{-}="prev"},
    @i;	@={"top","1245","1345","2345","123","12","13","23","45","bot"},
	@@{*{\mydot}},
	"12"+(-1.8,-1)*{{}^{1\!2}},
	"13"+(-1.5,-1)*{{}^{1\!3}},
	"23"+(2,-1)*{{}^{2\03}},
	"45"+(2,-1)*{{}^{4\05}},
	(77,0)="m",
	"m"+(16,-12)*{\mycaptionempty};
	"m"+(-14,0); "m"+(6,0)**\dir{-};
	"m"+(0,-6); "m"+(0,18)**\dir{-};
	"m"+(-11.4,-5.1); "m"+(3.6,17.4)**\dir{-};
	"m"+(0,0)*{\mydot};
	"m"+(0,12)*{\mydot};
	"m"+(-8,0)*{\mydot};
	"m"+(3.9,15.1)*{\scriptstyle 1\!2};
	"m"+(-1.5,-4.2)*{\scriptstyle 3\04};
	"m"+(-12.6,1.2)*{\scriptstyle 5\06};
	(16,-4)+"m"="m",
	(8,0)+"m"="bot",	
	(0,7)+"m"="12",
	(8,7)+"m"="34",
	(16,7)+"m"="56",
	(0,14)+"m"="1234",
	(8,14)+"m"="1256",
	(16,14)+"m"="3456",
	(8,21)+"m"="top",
    @i;	@={"1234","12","1256","56","3456","34","1234"},
	s0="prev" @@{;"prev";**@{-}="prev"},
    @i;	@={"1234","1256","3456"}, "top"; @@{**@{-}},
    @i;	@={"12","34","56"}, "bot"; @@{**@{-}},
    @i;	@={"top","1234","1256","3456","12","34","56","bot"},
	@@{*{\mydot}};
	"12"+(-1.8,-1)*{{}^{1\!2}},
	"34"+(2,-1)*{{}^{3\04}},
	"56"+(2,-1)*{{}^{5\06}},
\endxy
\end{figure}

Of special importance is the brick~$M^1_{1^r}$
that arises from the coordinate arrangement in~$\Bbbk^r$:
blow up $r$~points in $\P^{r-1}$ in general position,
then blow up the proper transforms of
all lines spanned by pairs of these points,
then blow up those of all planes spanned by triples,
and so on.
Thus $M^1_{1^r}$ is isomorphic to
the space~$\Pi_r$ that Kapranov called
the {\em permutahedral space} \cite[p.~105]{Kapranov:Chow}.
It~is the compact projective toric variety
whose encoding polytope is the permutahedron~$P_r$,
usually defined as the convex hull of the set of $r!$ points in~$`R^r$
with coordinates $(`s^{-1}(1),\dots,`s^{-1}(r))$,
for all $`s\in\Sym_r$.
This polytope can also be obtained from the standard $(r-1)$-simplex
by chopping off first all its vertices,
then all that remains of its edges,
then faces, and so on;
this corresponds to the sequence of blowups producing~$\Pi_r$.
In addition,
this variety is the closure of a principal toric orbit
in the complete flag variety
and it has been extensively studied from various perspectives
\cite{Atiyah,DolgachevLunts,GelfandSerganova,Procesi,Stanley:log-concave,%
Stembridge:Eulerian,Stembridge:reps}.

For each $m \> 1$,
the brick $M^m_{1^r}$ is a toric variety
because all strata of~$\Br^m_{1^r}$,
sitting in $\P^{rm-1}$,
are $(\kt)^{rm}$-invariant.
\end{examples}

\begin{prop}
\label{compound}
Every open compound brick has the structure of a bundle
\begin{equation}
\label{decompose}
\xymatrix @-3mm{
{}\0 M^1_{1^r} \ar[r] &
\0 M^m_{\lambda} \ar[r] &
\0 M^m_{\nu_1}\times\dots\times\0 M^m_{\nu_r},}
\end{equation}
where $`l$ is the integer partition $(`n_1,\dots,`n_r)$.
\end{prop}

\begin{proof}
The complement to $\Br^m_{`n_i}$ in $(\Am)^{`n_i}$ is $\Conf(\Am,`n_i+1)/\Am$,
the configuration space of $`n_i+1$ distinct labeled points in~$\Am$
modulo translations.
Since~$\0 M^m_`l$ is the complement in $\P(C^m_`l)$
to the projectivization of the arrangement
$\Br^m_`l = \Br^m_{`n_1} \times\dots\times \Br^m_{`n_r}$,
it follows that
\begin{equation}
\label{orbit}
\0 M^m_\lambda = \P(\0 C^m_\lambda), \quad\text{where}\quad
\0 C^m_\lambda = \prod_{i=1}^{r} \BBigfactor{\Conf(\Am,\nu_i+1)}{\Am},
\end{equation}
is the orbit space of the diagonal action of~$\kt$
on this product by dilations.

Separate actions of~$\kt$ on each factor
together give that of~$(\kt)^r$ on~$\0 C^m_`l$.
Its total orbit space is isomorphic to
the product of those coming from the factors,
which is $\0 M^m_{`n_1} \times\dots\times \0 M^m_{`n_r}$.
The orbit space $\0 M^m_`l$ maps into this product,
with fiber $(\kt)^r/\kt \simeq \0 M^1_{1^r}$.
\end{proof}

The next two propositions follow from
the general work of De Concini and Procesi \cite[pages~480--482]{DeCP},
but they can also be proved directly.

\begin{prop}
\label{affine}
\begin{enumerate}
\item
\label{affa}
The compactification $\Am\pd<n>$ is the product $\Am \times \L$,
where~$\L$ is the total space of a line bundle
over the simple brick~$M^m_{n-1}$.
\item
\label{affb}
The simple brick~$M^m_{n-1}$
is a compactification of\/ $\Conf(\Am,n)/\text{\rm Aff}$,
where $\text{\rm Aff}$ is the group of all affine transformations in~$\Am$.
\end{enumerate}
\end{prop}

\begin{proof}
(\ref{affa})
The direct factor~$\Am$ is the small diagonal $`D \subset (\Am)^n$.
The essential shape of the bottom partition of $[n]$ is the integer $n-1$,
thus by definition,
there is a map $\psi \colon M^m_{n-1} \to P = \P((\Am)^n/`D)$.
The bundle~$\L$ is the pullback by~$\psi$
of the tautologial line bundle over~$P$;
since~$\psi$ is an iterated blowup,
Lemma~\ref{blow-bundle}
(formulated below)
has to be used at each stage.

(\ref{affb})
Affine transformations identify any nondegenerate configuration in~$\Am$
with a degenerate one in which all $n$~points collide at $0$,
cancelling both the direct factor~$\Am$ and
the fiber of the line bundle~$\L^n$.
\end{proof}

\begin{lem}
\label{blow-bundle}
Let~$V$ be a smooth subvariety of a smooth algebraic variety~$W$, 
let $h\colon F\to W$ be a vector bundle over~$W$,
and~$E$ its restriction onto~$V$.
Then $\Bl_E F$ is a vector bundle over $\Bl_V W$
isomorphic to the pullback of~$F$ by the blowup projection.
\end{lem}

\begin{proof}
The normal bundle $\normal{E}{\!F}$ is
the pullback $h^*\normal{V}{W}$. 
\end{proof}

There are two differences between
the construction of $\Am\pd<n>$ and that of the bricks:
different arrangements to start with and projectivization;
both are minor enough that some basic facts about bricks
follow by the same arguments that apply to $\Am\pd<n>$.
In turn,
describing first the strata of the bricks
provides a quick way of doing the same for~$\Xn$.

\begin{prop}
\label{brick-strata}
For any integer partition~$`l$,
fix a partition $`p$ of $[n]$ of essential shape~$`l$
and an isomorphism $[`p,\top] \simeq L_{`l}$.
\begin{enumerate}
\item
\label{brick-divisors}
For each partition~$\1$,
$`p < \1 < \top$,
there is a divisor~$E^\1$ in~$M^m_`l$.
The union of these divisors is
the complement\/ $M^m_`l \smallsetminus \0 M^m_`l$,
and any set of them meets transversally.
\item
\label{brick-chains}
An intersection $E^{`p_1}\cap\dots\cap E^{`p_k}$
is nonempty exactly when the partitions form a chain.
Thus $M^m_`l$ is stratified by strata parametrized by
all chains in~$L_`l$ that include neither its bottom nor its top.
\item
\label{brick-products}
For any such chain $\{`p_1,\dots,`p_k\}$,
the corresponding stratum of\/~$M^m_`l$
is isomorphic to $M^m_{`l_0} \times\dots\times M^m_{`l_k}$,
where the integer partitions $`l_0,\dots,`l_k$
are determined by $L_{`l_i} \simeq [`p_i,`p_{i+1}]$,
with $`p_0 = `p$ and\/ $`p_{k+1} = \top$.
\end{enumerate}
\end{prop}

\begin{defn}
A smooth subvariety~$V$ of a smooth algebraic variety~$W$
will be called \newterm{straight}
if the normal bundle $N_V W$ is isomorphic to
a direct sum of copies of a single line bundle.
In this case,
the exceptional divisor of the blowup $\Bl_V W$ is a trivial bundle.
\end{defn}

\begin{lem}
\label{straight}
\begin{enumerate}
\item
\label{P-in-P}
For any two positive integers\/ $k$~and\/~$l$,
any linear subvariety\/ $\P^k$ of\/~$\P^{k+l+1}$ is straight.
{\rm (}\!Whence the term.{\rm )}
\item
\label{blow-straight}
Let\/ $Z$~and\/~$V$ be smooth subvarieties of a smooth algebraic variety~$W$,
such that either\/ $Z\cap V=\varnothing$ or $Z\subset V$.
If\/~$V$ is straight in~$W$,
then so is its proper transform~$\tilde{V}$ in $\tilde{W}=\Bl_Z W$.
\end{enumerate}
\end{lem}

\begin{proof}
Part (\ref{P-in-P}) follows directly from the definition.

(\ref{blow-straight})
Nothing to be done when $V$~and~$Z$ are disjoint.
When $Z\subset V$,
denote by~$E$ the exceptional divisor of $\tilde{W}$,
and by~$p$ the projection $\tilde{V}\to V$,
then
$$
\normal{\tilde{V}}{\tilde{W}} \simeq
p^*\normal{V}{W}\otimes\O(-E)\big|_{\tilde{V}},
$$
and the claim follows.
\end{proof}

\begin{proof}[Proof of Proposition~\ref{brick-strata}]

Similarly to Proposition~\ref{divisor},
the definition of~$M^m_`l$ implies parts
(\ref{brick-divisors})~and~(\ref{brick-chains}).

Part~(\ref{brick-products}) can be checked by induction on~$k$,
where the inductive step follows by applying the case $k=1$.
Thus,
it is enough to show that each divisor~$E^`p$
is isomorphic to $M^m_{`l_0} \times M^m_{`l_1}$,
where $L_{`l_0} \simeq [`p,\1]$ and $L_{`l_1} \simeq [\1,\top]$.
The argument is based on Lemmas \ref{blow-bundle}~and~\ref{straight}.

Every partition from $[`p,\top]$ belongs to one of the following six groups:
\begin{equation*}
\begin{align*}
\text{(i)}	\ &\ \{`p\},				&&&
\text{(ii)}	\ &\ \{`p' \suchthat `p<`p'<\1\},	&&\\
\text{(iii)}	\ &\ \{\1\},				&&&
\text{(iv)}	\ &\ \{`p' \suchthat \1<`p'<\top\},	&&\\
\text{(v)}	\ &\ \{\top\},				&&&
\text{(vi)}	\ &\ \text{incomparable with }~\1.
\end{align*}
\end{equation*}
The proof will be completed
by studying the impact of blowups corresponding to partitions in each group
on the stratum $`D^\1 / `D^`p$ of the arrangement $\Br^m_`l$.
Before the blowups,
the arrangements induced in $`D^\1 / `D^`p$ and $(\Am)^n / `D^\1$
are isomorphic respectively to $\Br^m_{`l_0}$ and $\Br^m_{`l_1}$.

First group,
first stage.
The exceptional divisor $\P(C^m_`l)$ of the first stage
has a straight subvariety $\P( `D^\1 / `D^`p) \simeq \P(C^m_{`l_0})$
with the projectivization of~$\Br^m_{`l_0}$ in it,
and with the arrangement $\Br^m_{`l_1}$ in each normal space to it
(Lemma~\ref{blow-bundle}).
Lemmas \ref{blow-bundle}~and~\ref{straight}
also apply at the subsequent stages,
pulling back arrangements inside normal spaces
and preserving the straightness of blowup centers.
Group~(vi) blowups are irrelevant for the divisor~$E^`p$ at all stages,
and no blowup corresponds to~$\top$.

Group~(ii) blowups turn $\P(C^m_{`l_0})$ into $M^m_{`l_0}$.
Then the group~(iii) blowup makes it into
a divisor isomorphic to $M^m_{`l_0} \times \P(C^m_{`l_1})$.
The second factor inherits the projectivization of~$\Br^m_{`l_1}$,
and blowups of the remaining group~(iv)
transform this divisor into $E^\1 \simeq M^m_{`l_0} \times M^m_{`l_1}$.
\end{proof}

\begin{lem}
\label{divisor-bundle}
Each divisor $D^`p$ of\/ $\Xn$ is isomorphic to 
a bundle over $X\pd<`r(`p)>$ with fiber~$M^m_{`l(`p)}$.
In addition,
this bundle is trivial if\/~$X=\Am$.
\end{lem}

\begin{proof}
Corollary~\ref{middleX} gives $Y^`p_{r-1} \simeq X\pd<r>$,
where $r =`r(`p)$.
By~Lemma~\ref{blow-bundle},
the arrangements~$\Br^m_`p$ transform isomorphically
from the normal spaces to $`D^`p$ in $\Xn$
into the normal spaces to $Y^`p_{r-1}$ in $Y^{\vphantom{`p}}_{r-1}$.
At the next stage,
$Y^`p_r$ is a bundle over $X\pd<r>$
with fibers isomorphic to $P = \P((\Am)^n/`D)$.
The relevant blowup centers of the subsequent stages are its subbundles;
their fibers form in every fiber of $Y^`p_r$ an arrangement
isomorphic to the projectivization of $\Br^m_`p$ in~$P$.
Thus in the end,
the fibers of $Y^`p_r$ transform into $M^m_{`l(`p)}$.

If in addition $X=\Am$,
a repeated application of Lemma~\ref{straight} shows that
$Y^`p_{r-1}$ is straight in $Y^`p$,
so $Y^`p_r = \Am\pd<r> \times P$ and the result follows.
\end{proof}

\begin{prop}
\label{tower}
Let\/~$`g = \{`p_1,\dots,`p_k\}$ be a chain of partitions of\/~$[n]$
and let\/ $\{`l_0,\dots,`l_k\}$ be
its associated sequence of integer partitions
(Section~\ref{sec-combinat}).
\begin{enumerate}
\item
The stratum~$S_`g$ of\/~$\Xn$
is a bundle over $X\pd<`l_0>$ with fiber isomorphic to
$M^m_{`l_1} \times\dots\times M^m_{`l_k}$.
\item
Consequently,
the complement in~$S_`g$ to the union of smaller strata,
the open stratum $\0 S_`g$,
is a bundle over $\Conf(X,`l_0)$ with fiber isomorphic to
$\0 M^m_{`l_1} \times\dots\times \0 M^m_{`l_k}$.
\end{enumerate}
\end{prop}

\begin{proof}
Put together Lemma~\ref{divisor-bundle} and Proposition~\ref{brick-strata}.
\end{proof}

By this proposition,
a~point in a stratum~$\0 S_`g$ 
is given by a~configuration of $r$~distinct points in~$X$
(where the collision occurs)
and a sequence consisting of one point in each open brick $\0 M^m_{`l_i}$.
Equations (\ref{decompose})~and~(\ref{orbit}) in Proposition~\ref{compound}
show that such points can be represented by suitable polyscreens:
points in each constituent open simple brick are Fulton--MacPherson screens,
and points in $\0 M^1_{1^r}$ are $r$-tuples of scale factors.
Thus,
points of~$\Xn$ indeed have the geometric description
explained in Section~\ref{sec-screens}.

\begin{prop}
\label{small-bricks}
The compound brick\/ $M^m_{1^r}$ has the structure of a bundle
\begin{equation}
\label{only-ones}
\xymatrix @-3mm{
\Pi_{r} \ar[r] &
M^m_{1^r} \ar[r] &
(M^m_1)^r. }
\end{equation}
\end{prop}

\begin{proof}
Fix a partition~$`p$ of $[2r]$ into two-element blocks
and let the nest $\S$ be the set $\{`b_1,\dots,`b_r\}$ of blocks of~$`p$.
The map $\vartheta_{2r} \colon \Am\pd<2r> \to \Am[2r]$
takes the divisor~$D^`p$ of $\Am\pd<2r>$
into the stratum $\Am(\S)$ of $\Am[2r]$.
The divisor is isomorphic to $\Am\pd<r> \times M^m_{1^r}$
by Lemma~\ref{divisor-bundle}
and the stratum is isomorphic to $\Am[r] \times (\P^{m-1})^r$.
Since $M^m_1 \simeq \P^{m-1}$,
it~follows that $M^m_{1^r}$ maps to $(M^m_1)^r$.

Tracing $\vartheta_{2r}^{-1}$ stage by stage,
first transform the factor $\Am[r]$ into $\Am\pd<r>$;
then at stage~$r$ blow up the proper transform of $\Am(\S)$,
turning the second factor into a $\P^{r-1}$-bundle over $(M^m_1)^r$.
Since $\Am[n] \smallsetminus \Conf(\Am,n)$
is a normal crossing divisor,
the divisors $D(`b_i)$ induce in each fiber $\P^{r-1}$
the projectivized coordinate hyperplane arrangement.
All of its strata are blown up at the subsequent stages,
turning $\P^{r-1}$ into $M^1_{1^r} \simeq \Pi_r$.
\end{proof}

\begin{figure}[t]
\label{levels-split}
\xy
	(-14,-17)*{\ };		(-14,9)*{\ };
	(80,-12)*{\mycaption{5in}{One level splits into two}},
	(0,0)="m"	*{\mydot};
     @i @={(0,3)+"m",(1.5,-3)+"m",(-1.5,-3)+"m"},
	"m"; @@{**@{-}},
	(6,0)+"m"="m"	*{\mydot};
     @i @={(0,3)+"m",(1.5,-3)+"m",(-1.5,-3)+"m"},
	"m"; @@{**@{-}},
	(6,0)+"m"="m"	*{\mydot};
     @i @={(0,3)+"m",(1.5,-3)+"m",(-1.5,-3)+"m"},
	"m"; @@{**@{-}},
	(6,0)+"m"="m"	*{\mydot};
     @i @={(0,3)+"m",(1.5,-3)+"m",(-1.5,-3)+"m"},
	"m"; @@{**@{-}},
	(6,0)+"m"="m"	*{\mydot};
     @i @={(0,3)+"m",(1.5,-3)+"m",(-1.5,-3)+"m"},
	"m"; @@{**@{-}},
	(6,0)+"m"="m"	*{\mydot};
     @i @={(0,3)+"m",(1.5,-3)+"m",(-1.5,-3)+"m"},
	"m"; @@{**@{-}},
	(6,0)+"m"="m"	*{\mydot};
     @i @={(0,3)+"m",(1.5,-3)+"m",(-1.5,-3)+"m"},
	"m"; @@{**@{-}};
	(-2,3.2);	 (38,3.2)**\dir{.};
	(-2,-3.2); (38,-3.2)**\dir{.};
	(18,5)*\txt{\scriptsize $\cdots$ upper levels $\cdots$};
	(18,-5)*\txt{\scriptsize $\cdots$ deeper levels $\cdots$};
%
	(60,0)="m",
     @i @={(0,4)+"m",(1.5,-4)+"m",(-1.5,-4)+"m"},
	"m"-(0,2); @@{**@{-}}; *{\mydot};
	(6,0)+"m"="m"
     @i @={(0,4)+"m",(1.5,-4)+"m",(-1.5,-4)+"m"},
	"m"+(0,2); @@{**@{-}}; *{\mydot};
	(6,0)+"m"="m"
     @i @={(0,4)+"m",(1.5,-4)+"m",(-1.5,-4)+"m"},
	"m"-(0,2); @@{**@{-}}; *{\mydot};
	(6,0)+"m"="m"
     @i @={(0,4)+"m",(1.5,-4)+"m",(-1.5,-4)+"m"},
	"m"+(0,2); @@{**@{-}}; *{\mydot};
	(6,0)+"m"="m"
     @i @={(0,4)+"m",(1.5,-4)+"m",(-1.5,-4)+"m"},
	"m"+(0,2); @@{**@{-}}; *{\mydot};
	(6,0)+"m"="m"
     @i @={(0,4)+"m",(1.5,-4)+"m",(-1.5,-4)+"m"},
	"m"-(0,2); @@{**@{-}}; *{\mydot};
	(6,0)+"m"="m"
     @i @={(0,4)+"m",(1.5,-4)+"m",(-1.5,-4)+"m"},
	"m"-(0,2); @@{**@{-}}; *{\mydot};
	"m"+(2,4.2);	"m"+(-38,4.2)**\dir{.};
	"m"+(2,-4.2);	"m"+(-38,-4.2)**\dir{.};
	"m"+(-18,6)*\txt{\scriptsize $\cdots$ upper levels $\cdots$};
	"m"+(-18,-6)*\txt{\scriptsize $\cdots$ deeper levels $\cdots$};
\endxy
\end{figure}

The fiber~$\Pi_r$ in Eq.~(\ref{only-ones}) stores scale factors;
points in its open part $\0 M^1_{1^r}$ are generic
and each is a part of one polyscreen.
The  divisor $E_r = \Pi_r \smallsetminus \0 M^1_{1^r}$
has components isomorphic to $\Pi_{s} \times \Pi_{r-s}$,
whose points represent degenerations
with $s$~scale factors tending to zero,
and therefore polyscreens that split into two:
$s$~screens form a~new level.
For example,
the left leveled tree in Figure~12 may degenerate into the right one,
corresponding to a divisor $\Pi_4 \times \Pi_3 \subset \Pi_7$.
The new levels may of course split further;
intersections of components of~$E_r$ give a stratification of~$\Pi_r$,
and each stratum is a product of
a number of smaller permutahedral varieties.
This corresponds to the well-known fact
that all faces of the permutahedron~$P_r$
are products of lower-dimensional permutahedra~\cite{BilleraSarangarajan}.

Other compound bricks,
that is,
$M^m_`l$ for those integer partitions~$`l$ that have parts greater than~$1$,
do not admit decompositions similar to Eq.~(\ref{only-ones}),
but each of them is a blowup of~$M^m_{1^r}$ for $r = |`l|$.
Let~$`L_r$ be the set of all partitions of an integer~$r$
partially ordered by refinement:
$(5,3) < (4,2,1,1)$ in~$`L_8$ because $5 = 4+1$ and $3 = 2+1$.
It turns out that
the set of bricks $\{M^m_`l \suchthat `l \in `L_r\}$
has a compatible (reverse) \inquotes{blowing-up} partial order.

\begin{prop}
\label{sub-and-super}
Suppose that\/ $`l,`l' \in `L_r$ and $`l <`l'$.
\begin{enumerate}
\item
\label{sublattice}
The lattice~$L_{`l'}$ contains a sublattice isomorphic to~$L_{`l}$.
\item
\label{subarrangement}
The subarrangement of~$\Br^m_{`l'}$ formed by the subspaces
that correspond to this sublattice is~$\Br^m_`l$,
up to coordinate change.
\item
\label{superbrick}
The brick~$M^m_{`l'}$ is an iterated blowup of\/~$M^m_{`l}$.
\end{enumerate}
\end{prop}

\begin{proof}
(\ref{sublattice})
It is enough to show this for $`l' = (r-1)$ and $`l = (s-1,r-s)$.
The required sublattice of~$L_{[r]}$
is generated by the union $[\1,\top] \cup [\2,\top]$,
where the only essential block of $\1$ ($\2$)
is $\{ k \suchthat k \< s \}$ (resp.~$\{ k \suchthat k \> s \}$).

(\ref{subarrangement})
It is enough to consider the same $`l$~and~$`l'$ as in~(\ref{sublattice})
and then write explicitly the equations for the large diagonals.

(\ref{superbrick})
The two lattices $L_`l \subset L_{`l'}$ determine
the sequences of blowups of $\P^{rm-1}$ creating $M^m_`l$ and $M^m_{`l'}$.
It suffices to show that
the blowups making~$M^m_{`l'}$ can be rearranged,
without changing the outcome (up to an isomorphism),
into a different sequence
so that an intermediate stage is $M^m_{`l}$.
This situation is quite similar to
the consideration of $\vartheta_n \colon \Xn\to X[n]$
in Section~\ref{sec-closure},
and similar is the solution.
\end{proof}

\section{Isotropy of the permutation action}
\label{sec-isotropy}

Assume that the ground field~$\Bbbk$ is of characteristic~$0$.
Reading carefully into \FM's proof of
the solvability of the isotropy subgroups of~$\Sym_n$ acting on~$X[n]$,
one soon realizes that
every point where the isotropy subgroup fails to be abelian
lies in a stratum whose encoding nest contains
a pair of disjoint subsets of~$[n]$.
Exactly these strata are blown up
by $\vartheta_n\colon\Xn\to X[n]$,
and this observation raises hopes that are not false.

\begin{thm}
\label{abelian}
If\/~$X$ is a smooth algebraic variety over a field\/~$\Bbbk$
of characteristic~$0$,
then all isotropy subgroups of\/~$\Sym_n$
acting on~$\Xn$ by permutations of labels
are abelian.
\end{thm}

\begin{proof}
First,
reduce to the case of all~$n$ points colliding at the same point in~$X$.
Suppose a collision~$\x$ occurs at $p_1,\dots,p_r\in X$.
If it could be studied near each~$p_i$ independently of the other points,
as for $X[n]$,
the isotropy subgroup would have been $G^{p_1}\times\dots\times G^{p_r}$,
where $G^{p_i}$ is the isotropy subgroup of the collision near~$p_i$.
It would have corresponded to~$r$ independent sequences of colored screens
and reduced the proof to the case of one collision point,
but this does not suit~$\Xn$.
Fortunately,
interdependencies among the corresponding levels in those $r$~sequences
only put more restrictions on a permutation aspiring to fix~$\x$.
It means that the isotropy subgroup will be a subgroup of the product above,
which still does the trick.

Pick a $k$-chain $`g\ni\bot$
and a coherent sequence of colored screens $\CS^j(\x)$ for~$`g$.
A~permutation $`s\in\Sym_n$ fixes~$\x\in\0 S_`g$
if and only if it fixes all $\CS^j(\x)$.
A~colored screen is fixed by~$`s$ exactly when
these two conditions are fulfilled:
\begin{enumerate}
\renewcommand\theenumi{\arabic{enumi}}
\renewcommand\labelenumi{(F\theenumi)}
\item
it is fixed modulo colors;
\item
any two points of the same color
go to two points of the same color,
not necessarily the original one.
\end{enumerate}

Let~$G$ be the isotropy subgroup at~$\x$.
A permutation $`s\in G$ satisfies~(F1) for $\CS^j(\x)$,
therefore it induces the scaling of $T_pX$ underlying $\CS^j(\x)$
by a scale factor $f_j(`s)\in\kt$.
The map $f_j\colon G\to\kt$ is a group homomorphism,
thus there is a group homomorphism $(f_1,\dots,f_k)= f\colon G\to (\kt)^k$,
and to show that it is injective suffices to complete the theorem.

Take $`s\in\ker f$,
then~$`s$ does not move points in any of the colored screens $\CS^j(\x)$,
$j=1,\dots,k$.
By coherence,
every color in $\CS^j(\x)$ is a point in $\CS^{j-1}(\x)$,
since both are but blocks of the partition $`p_j\in`g$.
Thus,
$`s$ cannot change colors either,
in any $\CS^j(\x)$ for $j=2,\dots,k$,
and colors in $\CS^1(\x)$ stay unchanged because there is only one such.
This shows that~$`s$ does not move anything at all,
and there is only one such permutation:
if $`s\neq\operatorname{id}$ and $`s(a)=b$,
then~$`s$ must induce nontrivial scaling on $\CS^l(\x)$,
where~$l$ is the maximal index~$j$ for which
$a$~and~$b$ are in the same block of $`p_j\in`g$.
\end{proof}

\begin{rem}
This version of the original proof
is one substantially simplified with a key idea due to Jean--Luc Brylinski.
\end{rem}


\end{document}